\documentclass[12pt, reqno]{amsart}
\usepackage{ifthen,amsfonts,amsmath,amssymb,epic,eepic}
\usepackage{epsfig,color}

\makeatother

\theoremstyle{definition}

\def\fnum{equation}
\newtheorem{Thm}[\fnum]{Theorem}
\newtheorem{Cor}[\fnum]{Corollary}
\newtheorem{Que}[\fnum]{Question}
\newtheorem{Lem}[\fnum]{Lemma}
\newtheorem{Con}[\fnum]{Conjecture}
\newtheorem{Def}[\fnum]{Definition}

\newtheorem{Pro}[\fnum]{Proposition}

\newcommand{\nn}{{\bf{n}}}

\newcommand{\diam}{{\text {diam}}}

\newcommand{\dist}{{\text {dist}}}

\def\ZZ{{\bold Z}}
\def\RR{{\bold R}}

\def\CC{{\bold C }}
\newcommand{\dv}{{\text {div}}}

\newcommand{\Cat}{{\text {Cat}}}

\newcommand{\K}{{\text{K}}}

\newcommand{\cB}{{\mathcal{B}}}

\newcommand{\cP}{{\mathcal{P}}}
\newcommand{\cS}{{\mathcal{S}}}

\newcommand{\rdel}{R_{\delta}}

\newcommand{\eqr}[1]{(\ref{#1})}

\newcommand{\cone}{{\bf{C}}}

\begin{document}

\title[The Calabi--Yau conjectures for embedded surfaces]
{The Calabi--Yau conjectures for embedded surfaces}

\author{Tobias H. Colding}%
\address{Courant Institute of Mathematical Sciences 251
Mercer Street\\ New York, NY 10012}
\author{William P. Minicozzi II}%
\address{Department of Mathematics\\
Johns Hopkins University\\
3400 N. Charles St.\\
Baltimore, MD 21218}
\thanks{The authors were partially supported by NSF Grants DMS
0104453 and DMS 0104187}


\email{colding@cims.nyu.edu and minicozz@math.jhu.edu}


\maketitle


\numberwithin{equation}{section}

\section{Introduction} \label{s:s0}

In this paper we will prove the Calabi-Yau conjectures for
embedded surfaces (i.e., surfaces without self-intersection).  In
fact, we will prove considerably more.  The heart of our argument
is very general and should apply to a variety of situations, as
will be more apparent once we describe the main steps of the proof
later in the introduction.

The Calabi-Yau conjectures about surfaces date back to the 1960s.
Much work has been done on them over the past four decades.  In
particular,    examples of Jorge-Xavier from 1980 and Nadirashvili
from 1996 showed that the immersed versions were false; we will
show here that for embedded surfaces, i.e., injective immersions,
they are in fact true.

Their original form  was given in 1965 in \cite{Ca} where E.
Calabi made the following two conjectures about minimal surfaces
(they were also promoted by S.S. Chern at the same time; see page
212 of \cite{Ch}):

\begin{Con} \label{con:1}
``Prove that a complete  minimal hypersurface in $\RR^n$ must be
unbounded.''
\end{Con}

Calabi continued:
``It is known that there are no compact minimal submanifolds of
$\RR^n$ (or of any simply connected complete Riemannian manifold
with sectional curvature $\leq 0$). A more ambitious conjecture
is'':

\begin{Con} \label{con:2}
``A complete [non-flat] minimal hypersurface in $\RR^n$ has an
unbounded projection in every $(n-2)$--dimensional flat
subspace.''
\end{Con}

These conjectures were revisited in S.T. Yau's 1982 problem list (see
 problem 91 in \cite{Ya1}) by which time the Jorge-Xavier paper had appeared:

\begin{Que}     \label{q:03}
``Is there any complete minimal surface in $\RR^3$ which is a subset
of the unit ball?

This was asked by Calabi, \cite{Ca}.  There is an example of a
complete [non-flat] minimally immersed surface between two parallel
planes due to L. Jorge and F. Xavier, \cite{JXa2}.  Calabi has
also shown that such an example exists in $\RR^4$. (One takes an
algebraic curve in a compact complex surface covered by the ball
and lifts it up.)''
\end{Que}

 The \underline{immersed}  versions of these conjectures turned out to be
 false.  As mentioned above,  Jorge and Xavier, \cite{JXa2},  constructed
non-flat minimal immersions contained between two parallel planes
in 1980, giving a counter-example to the immersed version of the
more ambitious Conjecture \ref{con:2}; see also \cite{RoT}.
 Another significant development came in 1996,  when
 N. Nadirashvili, \cite{Na1}, constructed a complete
immersion of a minimal disk into the unit ball in $\RR^3$, showing
that Conjecture \ref{con:1} also failed for immersed surfaces; see
  \cite{MaMo1}, \cite{LMaMo1}, \cite{LMaMo2},
 for other topological types than disks.

\vskip2mm
 The conjectures were again revisited in Yau's 2000
millenium lecture (see page 360 in \cite{Ya2}) where Yau stated:

\begin{Que}     \label{q:04}
``It is known \cite{Na1} that there are complete minimal surfaces
properly immersed into the [open] ball.  What is the geometry of
these surfaces?  Can they be embedded?...''
\end{Que}

As mentioned in the very beginning of the paper, we will in fact
show considerably more  than Calabi's conjectures.  This is in
part because the conjectures  are closely related to properness.
Recall that an immersed surface in an open subset $\Omega$ of
Euclidean space $\RR^3$ (where $\Omega$ is  all of $\RR^3$ unless
stated otherwise) is {\it proper} if the pre-image of any compact
subset of $\Omega$ is compact in the surface. A well-known
question generalizing Calabi's first conjecture asks when is a
complete \underline{embedded} minimal surface proper? (See for
instance question $4$ in \cite{MeP}, or the ``Properness
Conjecture'', conjecture $5$, in \cite{Me}, or question $5$ in
\cite{CM7}.)

Our main result  is  a chord arc bound{\footnote{A chord arc bound
is a bound from above and below for the ratio of intrinsic to
extrinsic distances.}} for intrinsic balls   that implies
 properness.
Obviously, intrinsic distances are larger than extrinsic
distances, so the significance of a chord arc bound is the reverse
inequality, i.e., a bound on intrinsic distances from above by
extrinsic distances. This is accomplished in the next theorem:

\begin{Thm}     \label{t:1}
There exists a  constant  $C > 0$  so that if $\Sigma \subset
\RR^3$ is an embedded minimal disk, $\cB_{2R}=\cB_{2R}(0)$ is an
intrinsic ball{\footnote{Intrinsic balls will be denoted with
script capital ``b'' like $\cB_r(x)$ whereas extrinsic balls will
be denoted by an ordinary capital ``b'' like $B_r(x)$.}}
 in $\Sigma \setminus
\partial \Sigma$ of radius $2R$, and if
$ \sup_{\cB_{r_0}}|A|^2>r_0^{-2}$ where $R>  r_0$,
 then for $x \in \cB_R$
 \begin{equation}   \label{e:t1}
    C \, \dist_{\Sigma}(x,0)<   |x| + r_0   \, .
 \end{equation}
\end{Thm}

The assumption of a lower bound for the supremum of the sum of the
squares of the principal curvatures, i.e., $
\sup_{\cB_{r_0}}|A|^2>r_0^{-2}$, in the theorem is a necessary
normalization for a chord arc bound.  This can easily be seen by
rescaling and translating the helicoid.  Equivalently this
normalization can be expressed in terms of the curvature, since by
the Gauss equation $-\frac{1}{2}|A|^2$ is equal to the curvature
of the minimal surface.

 Properness of a complete embedded minimal disk   is an immediate
consequence of Theorem \ref{t:1}.  Namely, by \eqr{e:t1}, as
intrinsic distances go to infinity, so do extrinsic distances.
Precisely, if $\Sigma$ is flat, and hence a plane, then obviously
$\Sigma$ is proper and if it is non-flat, then
$\sup_{\cB_{r_0}}|A|^2>r_0^{-2}$ for some $r_0>0$ and hence
$\Sigma$ is proper by \eqr{e:t1}.  In sum, we get the following
corollary:

\begin{Cor}     \label{c:2nd}
A complete embedded minimal disk in $\RR^3$ must be proper.
\end{Cor}

Corollary \ref{c:2nd} in turn implies that the first of
  Calabi's conjectures is true for
{\underline{embedded}} minimal disks. In particular,
Nadirashvili's examples cannot be embedded.
 We also get from it an answer to Yau's questions
 (Question \ref{q:03} and Question \ref{q:04}).

Another immediate consequence of Theorem \ref{t:1} together with
the one-sided curvature estimate of \cite{CM6} (i.e., theorem
$0.2$ in \cite{CM6}) is the following  version of that estimate
for intrinsic balls; see question $3$ in \cite{CM7} where this was
conjectured:

\begin{Cor}  \label{t:one-sided}
There exists $\epsilon>0$, so that if
\begin{equation}
    \Sigma \subset \{x_3>0\} \subset \RR^3
\end{equation}
 is an embedded
minimal disk with intrinsic ball $\cB_{2R} (x) \subset \Sigma
\setminus
\partial \Sigma$ and $|x|<\epsilon\,R$,  then
\begin{equation}        \label{e:graph}
\sup_{ \cB_{R}(x) } |A_{\Sigma}|^2 \leq R^{-2} \, .
\end{equation}
\end{Cor}

 As a corollary of this intrinsic one-sided curvature estimate we get that
the second, and ``more ambitious'', of Calabi's conjectures is
also true for {\underline{embedded}} minimal disks.  In
particular,
 Jorge-Xavier's examples cannot be embedded. Namely, letting  $R \to \infty$
 in Corollary  \ref{t:one-sided} gives the following halfspace
 theorem:

\begin{Cor}     \label{c:1}
The plane is the only complete embedded minimal disk in $\RR^3$
 in a halfspace.
\end{Cor}

In the final section, we will see that our results for disks imply
both of Calabi's conjectures and properness also for embedded
surfaces with finite topology. Recall that a surface $\Sigma$ is
said to have finite topology if it is homeomorphic to a closed
Riemann surface with a finite set of points removed or
``punctures''.  Each puncture corresponds to an end of $\Sigma$.

The following generalization of the halfspace theorem gives
Calabi's second, ``more ambitious'', conjecture for embedded
surfaces with finite topology:

\begin{Cor}     \label{c:finite}
The plane is the only complete embedded minimal surface with
finite topology
 in a halfspace of
 $\RR^3$.
\end{Cor}

Likewise, we get the properness of embedded surfaces with finite
topology:

\begin{Cor}     \label{c:3rd}
A complete embedded minimal surface with finite topology in
$\RR^3$ must be proper.
\end{Cor}

Most of the classical theorems on minimal surfaces assume
properness, or something which implies properness (such as finite
total curvature).  In particular, this assumption can now be
removed from these theorems.

\vskip2mm
Before we recall in more detail some of the earlier work on
these conjectures we will try to give the reader an idea of why these
kind of properness results should hold.

The proof that complete
embedded minimal disks are proper, i.e., Corollary \ref{c:2nd},
consists roughly of the
following three main steps:

\begin{enumerate}
\item
Show that if the surface is compact in a ball, then in this
ball we have good chord arc bounds.
\item
Show that if each component of the intersection of each ball of a
certain size is compact
(so that by the first step we have good estimates),
then each intersection with double the Euclidean balls is also compact.
Initially possibly with a much worse constant but then by the
first step with a good constant.
\item
Iterate the above two steps.
\end{enumerate}

\vskip2mm Step 1 above relies on our earlier results (see
\cite{CM3}--\cite{CM6}; see also \cite{CM9} for a survey) about
\underline{properly} embedded minimal disks.  We will come back to
this in the main body of this paper and instead here outline the
proof of step 2 assuming step 1.

\vskip2mm
 Suppose therefore that all intersections of the given disk with all Euclidean balls
of radius $r$ are compact and have good chord arc bounds.
We will show the same for all Euclidean balls of radius $2r$.

If not; then there are two points $x$, $y\in B_{2r}\cap \Sigma$
in the same connected component of $B_{2r}\cap \Sigma$
but with $\text{dist}_{\Sigma}(x,y)\geq C\, r$ for some
large constant $C$. Let $\gamma$ be an intrinsic geodesic in
$B_{2r}\cap \Sigma$
connecting $x$ and $y$.
 By dividing $\gamma$ into segments, we conclude that there must be a
pair of points $x_0$ and $y_0$ on $\gamma$ in $B_{2r}$ which are
intrinsically far apart yet extrinsically close.
 We will start at these two points and build out showing that $x_0$
and $y_0$ could not connect in $B_{2r}\cap \Sigma$.  This will be
the desired contradiction.

 By the assumption, each component of $B_r(x_0)\cap \Sigma$
is compact and by step
1 has good chord arc bounds; hence $x_0$ and $y_0$ must lie in
different components. Thus we have two compact components
of $B_{r}(x_0)\cap \Sigma$ which are
extrinsically close near the center.
 Earlier results (the one-sided curvature estimate of
\cite{CM6}; see theorem 0.2 there) show that half of each of these
two components must have curvature bounds.  Since this bound for
the curvature is in terms of the size of the relevant balls, then
it follows that a fixed fraction of these components must be
almost flat - again relative to its size.  In fact, it follows now
easily that these two almost flat regions contains intrinsic balls
centered at $x_0$ and $y_0$ and with radii a fixed fraction of
$r$.
  We can therefore go to the boundary of these almost flat intrinsic
balls and find two points
$x_1$ and $y_1$; one point in each intrinsic ball which are
extrinsically close yet intrinsically far apart.

Repeat the argument with $x_1$ and $y_1$ in place of $x_0$ and
$y_0$ to get points $x_2$ and $y_2$.
 Iterating gives large regions in the surface
centered at $x_0$ and $y_0$ with a priori
curvature bounds.
 Once we have a priori curvature bounds then
improvements involving stability show that even these large
regions are almost flat and thus could not combine in $B_{2r}$.
 This is the desired contradiction and hence completes the outline of
step 2 above of the proof that embedded minimal disks are proper.

\vskip6mm
 It is clear from the definition of proper that a proper
minimal surface in $\RR^3$ must be unbounded, so the examples of
Nadirashvili are not proper.  Much less obvious is that the plane
is the only complete {\underline{proper}} immersed minimal surface
in a halfspace.  This is however a consequence of the strong
halfspace theorem of D. Hoffman and W. Meeks, \cite{HoMe}, and
implies that also the examples of  Jorge-Xavier are not proper.

There has been extensive  work on both properness (as in Corollary
\ref{c:2nd}) and the halfspace property (as in Corollary
\ref{c:1})   assuming various {\underline{curvature}}
{\underline{bounds}}. Jorge and Xavier, \cite{JXa1} and
\cite{JXa2}, showed that there cannot exist a complete immersed
minimal surface with {\underline{bounded}} {\underline{curvature}}
in $\cap_i \{ x_i > 0 \}$; later Xavier proved that the plane is
the only such surface in a halfspace, \cite{Xa}. Recently, G.P.
Bessa, Jorge and G. Oliveira-Filho, \cite{BJO}, and H. Rosenberg,
\cite{Ro}, have shown that if complete  embedded minimal surface
has bounded curvature, then it must be proper. This properness was
extended to embedded minimal surfaces with locally bounded
curvature and finite topology by Meeks and Rosenberg in
\cite{MeRo}; finite topology was subsequently replaced by finite
genus in \cite{MePRs} by Meeks, J. Perez and A. Ros.

Inspired by Nadirashvili's examples, F. Martin and S. Morales
 constructed  in \cite{MaMo2} a complete bounded minimal
immersion which is proper in the (open) unit ball.  That is, the
preimages of compact subsets of the (open) unit ball are compact
in the surface and the image of the surface accumulates on the
boundary of the unit ball. They extended this  in
 \cite{MaMo3} to show that any convex, possibly noncompact or
nonsmooth, region of $\RR^3$ admits a proper complete minimal
immersion of the unit disk; cf.  \cite{Na2}.

Finally, we note that Calabi and P. Jones, \cite{Jo}, have
constructed bounded complete holomorphic (and hence minimal)
\underline{embeddings} in higher codimension.  Jones' example is a
graph and he used purely analytic
methods (including the Fefferman-Stein duality theorem between $H^1$ and BMO)
while, as mentioned in Question \ref{q:03}, Calabi's
approach was algebraic: Calabi considered the lift of an algebraic
curve in a complex surface covered by the unit ball.

\vskip2mm Throughout this paper, we
 let $x_1 , x_2 , x_3$ be the standard coordinates on
$\RR^3$. For $y \in  \Sigma \subset \RR^3$ and $s > 0$, the
extrinsic and intrinsic balls are $B_s(y)$ and $\cB_s(y)$,
respectively,   and  $\dist_{\Sigma} (\cdot , \cdot)$ is the
intrinsic distance in $\Sigma$.  We will use $\Sigma_{y,s}$ to
denote the component of $B_{s}(y) \cap \Sigma$ containing $y$; see
Figure \ref{f:cy1}. The two-dimensional disk $B_s(0) \cap \{ x_3 =
0 \}$ will be denoted by $D_s$.  The sectional curvature of a
smooth surface $\Sigma\subset \RR^3$ is $\K_{\Sigma}$  and
$A_{\Sigma}$ will be its second fundamental form. When $\Sigma$ is
oriented, $\nn_{\Sigma}$ is the unit normal.

\begin{figure}[htbp]
    \setlength{\captionindent}{20pt}
    \centering\input{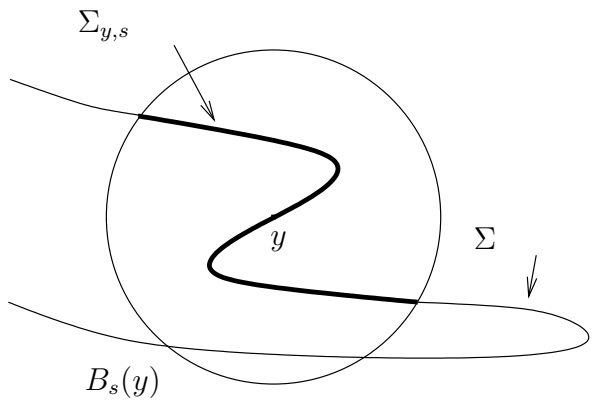}
    \caption{$\Sigma_{y,s}$ (in bold) denotes the component of $B_{s}(y) \cap
\Sigma$ containing $y$.}
    \label{f:cy1}
\end{figure}

We will use freely that each component of the intersection of a
minimal disk with an extrinsic ball is also a disk (see, e.g.,
appendix C in \cite{CM6}).  This follows easily from the maximum
principle since $|x|^2$ is subharmonic on a minimal surface.

\vskip2mm In \cite{CM9}, the results of this paper as well as
\cite{CM3}--\cite{CM6} are surveyed.

\section{Theorem \ref{t:1} and estimates for intrinsic balls}

The main result of this paper (Theorem \ref{t:1})  will follow by
combining the next proposition with a result from \cite{CM6}. This
next proposition gives a weak chord arc bound for an embedded
minimal disk but, unlike Theorem \ref{t:1}, only for one component
of a smaller extrinsic ball.  The result from \cite{CM6} will then
be used to show that there is in fact only one component, giving
the theorem.

\begin{Pro}     \label{t:2}
There exists $\delta_1 > 0$ so that if $\Sigma \subset \RR^3$ is
an embedded minimal disk, then for all intrinsic balls
 $\cB_R(x)$ in $ \Sigma \setminus
\partial \Sigma$,  the component $\Sigma_{x,\delta_1 \, R}$ of
$B_{\delta_1 \, R}(x) \cap \Sigma$ containing $x$ satisfies
 \begin{equation}   \label{e:t2}
    \Sigma_{x,\delta_1 \, R} \subset \cB_{R/2}(x) \, .
 \end{equation}
\end{Pro}

The result that we need from \cite{CM6} to show Theorem \ref{t:1}
is a consequence  of the one-sided curvature estimate of
\cite{CM6}; it is corollary $0.4$ in \cite{CM6}. This corollary
says that if two disjoint embedded minimal disks with boundary in
the boundary of a ball both come close to the center, then each
has an interior curvature estimate. Precisely, this is the
following result:

\begin{Cor}      \label{c:barrier}
 \cite{CM6} There exist constants $c > 1$ and $\epsilon >0$ so that the
following holds:

\noindent Let $\Sigma_1$ and $ \Sigma_2$ be disjoint embedded
minimal surfaces in  $B_{cR} \subset \RR^3$ with $\partial
\Sigma_i \subset
\partial B_{cR}$ and $B_{\epsilon \, R } \cap \Sigma_i \ne
\emptyset$. If $\Sigma_1 $ is a disk,
 then for all components $\Sigma_1'$ of
$B_{R} \cap \Sigma_1$ which intersect $B_{\epsilon \, R}$
\begin{equation}        \label{e:onece}
    \sup_{\Sigma_1'}   |A|^2
        \leq  R^{-2}  \, .
\end{equation}
\end{Cor}

Using this corollary, we can now prove Theorem \ref{t:1} assuming
Proposition \ref{t:2}, whose proof will fill up the next two
sections.

\begin{proof}
(of Theorem \ref{t:1} using Corollary \ref{c:barrier} and assuming
Proposition \ref{t:2}). Let $c > 1$ and $\epsilon > 0$ be given by
Corollary \ref{c:barrier} and $\delta_1 > 0$ by Proposition
\ref{t:2}.

Let $x \in \cB_R(0)$ be a fixed but arbitrary point and let
$\Sigma_0$ and $\Sigma_x$ be the components of
\begin{equation}
    B_{ \frac{c \, (|x| + r_0)}{\epsilon} } \cap \Sigma
\end{equation}
containing $0$ and $x$, respectively.  Here $r_0$ is given by the
curvature assumption in the statement of the theorem.  We will
divide into two cases depending on whether or not we have the
following inequality
\begin{equation}    \label{e:Rbig}
    \frac{2\, c \, (|x| + r_0)}{\delta_1 \, \epsilon} \leq R \, .
\end{equation}

If \eqr{e:Rbig} holds, then Proposition \ref{t:2} (with radius
equal to $\frac{2\, c \, (|x| + r_0)}{\delta_1 \, \epsilon}$)
implies that
\begin{equation}        \label{e:t2a}
    \Sigma_0  \subset \cB_{ \frac{c \, (|x| + r_0)}{\delta_1 \, \epsilon} }
    (0)
\end{equation}
and also, since $B_{ \frac{c \, (|x| + r_0)}{\epsilon} } \subset
B_{ \frac{2\, c \, (|x| + r_0)}{\epsilon} }(x)$ by the triangle
inequality,
\begin{equation}        \label{e:t2b}
    \Sigma_x  \subset \cB_{ \frac{c \, (|x| + r_0)}{\delta_1 \, \epsilon}
     }
    (x) \, .
\end{equation}
On the other hand, by definition, the embedded minimal disks
$\Sigma_0$ and $\Sigma_x$ are contained in $B_{ \frac{c \, (|x| +
r_0)}{\epsilon} }$.  Since $0$ and $x$ are in the smaller
extrinsic ball $B_{ c \, (|x| + r_0) }$, then both $\Sigma_0$ and
$\Sigma_x$ intersect $B_{ c \, (|x| + r_0) }$.
  Furthermore,
\eqr{e:t2a} and \eqr{e:t2b} imply that $\Sigma_0$ and $\Sigma_x$
are both compact and have boundary in $\partial B_{ \frac{c \,
(|x| + r_0)}{\epsilon} }$.  However, it follows from Corollary
\ref{c:barrier} and the lower curvature bound (i.e.,
$\sup_{\cB_{r_0}}|A|^2>r_0^{-2}$)   that there can only be one
 component with all of these properties.
 Hence, we have $\Sigma_0 = \Sigma_x$ so that
\begin{equation}    \label{e:x}
    \Sigma_x \subset \cB_{ \frac{c \, (|x| + r_0)}{ \delta_1 \, \epsilon} }
    (0) \, ,
\end{equation}
giving the claim \eqr{e:t1}.

 In the
remaining case, where \eqr{e:Rbig} does not hold, the claim
\eqr{e:t1} follows trivially.
\end{proof}

\vskip2mm
Before discussing the proof of Proposition \ref{t:2}, we
conclude this section by noting some additional applications of
Theorem \ref{t:1}.
 As alluded to in the introduction, then an immediate
consequence of Theorem \ref{t:1} is that we get   intrinsic
versions of all of the results of \cite{CM6}.  For instance we get
the following:

\begin{Thm}     \label{t:maincm6}
Intrinsic balls in embedded minimal disks are part of  properly
embedded double spiral staircases.  Moreover, a sequence of such
disks with curvature blowing up converges to a lamination.
\end{Thm}

For a precise statement of Theorem \ref{t:maincm6}, see theorem
$0.1$ of \cite{CM6}, with  intrinsic balls instead of extrinsic
balls.

A double spiral staircase consists of two multi-valued graphs (or
 spiral staircases) spiralling together around a common axis, without intersecting, so that
 the the flights of stairs alternate between the two staircases.
  Intuitively, an (embedded) multi--valued
graph is a surface such that over each point of the annulus, the
surface consists of $N$ graphs; the actual definition is recalled
in Appendix \ref{s:a1}.

\section{Chord arc properties of {\underline{properly}} embedded minimal disks}
\label{s:s2}

The proof of Proposition \ref{t:2} will be divided into several
steps over the next two sections.  The first step is to prove the
special case where we assume in addition that $\Sigma$ is compact
and has boundary in the boundary of an extrinsic ball.  The
advantage of this assumption is that the results of
\cite{CM3}-\cite{CM6} can   be applied directly.

\subsection{Properly embedded
disks}  The next proposition gives a   weak chord arc bound for a
compact  embedded minimal disk with boundary in the boundary of a
ball.  The fact that this  bound is otherwise independent of
$\Sigma$ will be crucial later when remove these assumptions.

\begin{Pro}     \label{t:c-a}
Let $\Sigma \subset   \RR^3$ be a {\underline{compact}} embedded
minimal disk.  There exists a constant $\delta_2 > 0$ independent
of $\Sigma$  such that if $x \in \Sigma$ and $\Sigma \subset B_R
(x)$ with  $\partial \Sigma \subset
\partial B_R(x)$, then the component
$\Sigma_{x,\delta_2 \, R}$ of $B_{\delta_2 \, R} (x) \cap \Sigma$
containing $x$ satisfies
 \begin{equation}   \label{e:c-a}
    \Sigma_{x,\delta_2 \, R} \subset \cB_{ \frac{R}{2} }(x) \, .
 \end{equation}
\end{Pro}

The key ingredient in the proof of  Proposition \ref{t:c-a} is
an effective version of the first main theorem in \cite{CM6}.
Before we can state this effective version, we need to
recall two definitions from \cite{CM6}.

First, given a constant $\delta > 0$ and a point $z \in \RR^3$,
then we denote by $\cone_{\delta}(z)$ the (convex) cone with
vertex $z$, cone angle $(\pi/2 - \arctan \delta)$, and axis
parallel to the $x_3$-axis. That is,
\begin{equation}
\cone_{\delta}(z)=\{x\in \RR^3\,|\,(x_3-z_3)^2 \geq
\delta^2\,((x_1-z_1)^2+(x_2-z_2)^2) \} \, .
\end{equation}

Second, recall from \cite{CM6} that,  roughly speaking, a blow up
pair $(y,s)$ consists of a point $y$ where the curvature is almost
maximal in a (extrinsic) ball of radius roughly $s$.   To be
precise,
 fix a constant $C_1$, then a point $y$ and scale $s > 0$ is a {\it blow up
pair} $(y,s)$  if
\begin{equation}        \label{e:int222}
    \sup_{B_{C_1 \, s}(y)\cap \Sigma} |A|^2 \leq 4 \, s^{-2} = 4 \, |A|^2 (y) \,
    .
\end{equation}
The constant $C_1$ will be given by theorem $0.7$ in \cite{CM6}
that gives the existence of a multi-valued graph starting on the
scale $s$.

\vskip2mm
 We are now ready to state a local version of
the first main theorem in \cite{CM6}. This is Lemma \ref{l:t0.1}
below and
 shows that a compact embedded minimal
disk, with boundary in the boundary of an extrinsic ball, is part
of a double spiral staircase.  In particular, it consists of two
multi-valued graphs spiralling together away from    a collection
of balls whose centers  lie along a Lipschitz curve transverse to
the graphs.  (The centers $y_i$ will be ordered by height around a
``middle point'' $y_0$;  negative values of $i$ should be thought
of as points below $y_0$.)

\begin{Lem} \label{l:t0.1}
Let $\Sigma \subset \RR^3$   be a compact embedded minimal disk.
There exist constants $c_{in}$, $c_{out}$, $c_{dist}$, $c_{max}$,
and $\delta
> 0$ independent of $\Sigma$ so that if
 $\Sigma \subset B_{R}$ with $\partial \Sigma \subset
\partial B_{R}$ and
\begin{equation}    \label{e:cmax}
    \sup_{B_{  R /c_{max} } \cap \Sigma}|A|^2 \geq c_{max}^2 \, R^{-2} \, ,
\end{equation}
 then there is  a collection of blow up pairs $\{ (y_i , s_i)
 \}_{i}$ with  $y_0 \in B_{R/(4c_{out})}$.  In addition,
after a rotation of $\RR^3$, we have that:
 \begin{enumerate}
 \item[(0)] For every $i$, we have $B_{C_1 \, s_i}(y_i)
\subset B_{6 R/ c_{out}}$.
 \item The extrinsic balls  $B_{s_i}(y_i)$ are disjoint and the points $\{
y_i \}$ lie in the intersections of the cones
\begin{equation}
    \cup_i \{ y_i \} \subset \cap_i \cone_{\delta}(y_i) \, .
\end{equation}
\item
  The points $y_i$ ``string together'' starting at $y_0$: For each $i>0$,
we have $y_{i}  \in B_{c_{in} \,
  s_i}(y_{i-1})$; for each $i< 0$, we have $y_{i}  \in B_{c_{in} \,
  s_i}(y_{i+1})$.
\item
    The $y_i$'s go from top to bottom, i.e., there is a curve $\tilde{\cS}
    \subset B_{R/c_{out}} \cap \cup_i B_{c_{in} \, s_i}(y_i)$ with
\begin{equation}
    \inf_{\tilde{\cS}} x_3 \leq - \frac{\delta\,R}{2 \, c_{out}} <
        \frac{\delta\, R}{2 \, c_{out}}  \leq \sup_{\tilde{\cS}}
    x_3 \, .
\end{equation}
\item ``Graphical away from balls'':
$B_{R/c_{out}} \cap \Sigma \setminus \cup_i B_{c_{in} \, s_i}(y_i)
$ consists of exactly two multi-valued graphs
 (which spiral together) with gradient $\leq \delta / 2$.
 \item
  ``Chord arc'': For each $i$, we have $B_{c_{in} \, s_i}(y_i) \cap
\Sigma \subset \cB_{c_{dist} \,
  s_i}(y_i)$.
\end{enumerate}
\end{Lem}

\begin{figure}[htbp]
    \setlength{\captionindent}{20pt}
    \centering\input{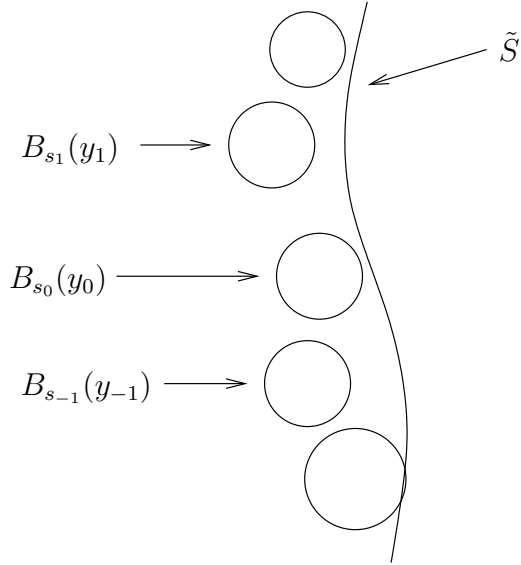}
    \caption{The balls $B_{s_i}(y_i)$ in the statement of Lemma \ref{l:t0.1}
are disjoint, yet consecutive balls are not too far apart; cf (2).
 In particular, the ratio of the radii of consecutive balls is bounded.}
    \label{f:cy1}
\end{figure}

Note that (1)--(3) are the effective version of the fact that the
singular set $\cS$ in \cite{CM6} is a Lipschitz graph over the
$x_3$-axis.    Property (4) says that the surface is graphical
away from the balls $B_{c_{in} \, s_i}(y_i) $.  Finally, (5) is a
chord arc property showing that the extrinsic balls $B_{c_{in} \,
s_i}(y_i) $ are contained in intrinsic balls $\cB_{c_{dist} \,
s_i}(y_i) $.

The proof of Lemma \ref{l:t0.1} is essentially contained in
\cite{CM6} but was not made explicit there.  We will describe
where to find properties (0)-(5) in \cite{CM6}, as well as the
necessary modifications, over the next three subsections.  The
reader who wishes to take these six properties (0)-(5) for granted
should jump to subsection \ref{ss:25}.

\subsection{Results from \cite{CM6}}

We will first recall a few of the results from \cite{CM6} that we
will use. The first of these, theorem 0.7 in \cite{CM6}, gives the
existence of multi-valued graphs nearby a blow up pair; cf.
\eqr{e:int222}. The precise statement is the following:

\begin{Lem} \label{t:mv} \cite{CM6}
Given $N\in \ZZ_+$ and $\epsilon > 0$, there exist  $C_1$ and $C_2
>0$ so that the following holds:

\noindent Let $\Sigma \subset \RR^3$ be an embedded minimal disk
with $0\in \Sigma \subset B_{R}$ and $\partial \Sigma\subset
\partial B_{R}$. If $(0,s)$ with $0<s<R/C_1$
is a blow up pair (i.e., satisfies \eqr{e:int222} with $y=0$ and
this $C_1$),
  then there exists (after a rotation of $\RR^3$) an $N$-valued
graph
\begin{equation}
    \Sigma_g \subset \Sigma\cap \{ x_3^2 \leq
    \epsilon^2 \, (x_1^2 + x_2^2) \}
\end{equation}
 over $D_{R/C_2}
\setminus D_{C_1 s}$ with gradient $\leq \epsilon$.
\end{Lem}

The second result that we will need to recall is the existence of
blow up pairs nearby a given blow up pair. This will be used to
show that the points $y_i$ string together.   This was a key
ingredient in the proofs of both main theorems in \cite{CM6} and
is recorded in proposition I.0.11 there (it was proven in
corollary III.3.5 in \cite{CM5}). For clarity, we restate this
next and give an elementary proof using \cite{CM6}.  Note,
however, that we could not have used this elementary proof in
\cite{CM5} since \cite{CM6} relies on \cite{CM5}.

\begin{Lem}     \label{l:nearby}    \cite{CM5}
Let $N$, $\epsilon$, $C_1$, and $C_2$ be as in  lemma \ref{t:mv}.
Then there exists a constant $C_5
> 4 \, C_1$ so that if
\begin{enumerate}
\item[(a)]
$\Sigma \subset \RR^3$ is an embedded minimal disk with $\Sigma
\subset B_{C_5 \, s}(y)$ and $\partial \Sigma \subset
\partial B_{C_5 \, s}(y)$,
\item[(b)]
$(y,s)$ is a blow up pair,
\end{enumerate}
then we get two blow up pairs $(y_+ , s_+)$ above $y$ and $(y_- ,
s_-)$ below $y$ with
\begin{equation}    \label{e:bupin}
    B_{C_1 \, s_{\pm}} (y_{\pm}) \subset B_{C_5 \, s}(y) \setminus
    B_{C_1 \, s}(y) \, .
\end{equation}
\end{Lem}

\begin{proof}
After rescaling and translating  $\Sigma$, we can assume
 that $y =0$ and $s =1$.  We will find the blow up pair $(y_+ , s_+)$ above $y$
 (the other case is identical).  Let $\Sigma^+$ denote the portion of $\Sigma$ above $0$
 (i.e., above the multi-valued graph corresponding to this blow up pair).

   It is easy to see by a simple blow up argument
 (lemma $5.1$ in \cite{CM4}) that it suffices to show that
\begin{equation}   \label{e:c5}
    \sup_{z \in B_{C_5/2} \cap \Sigma^+ \setminus B_{4 \, C_1 } } \,
    |z|^{-2} \,  |A|^2 (z)  \geq  4\, C_1  \, .
\end{equation}
We will argue by contradiction; suppose therefore that $\Sigma_i$
 is a sequence of embedded minimal disks
  satisfying (a) and (b) with $y=0$, $s=1$, and $C_5 = i$
 but so that \eqr{e:c5} fails for every $i$.

 Rescaling the $\Sigma_i$'s by a factor of $\sqrt{i}$, we get a
 new sequence $\tilde{\Sigma}_i$ with
$\tilde{\Sigma}_i \subset B_{\sqrt{i}}$   and
     $\partial \tilde{\Sigma}_i \subset \partial B_{\sqrt{i}}$ and
     so that
 $|A|^2(0) \to \infty$.
   Hence,
 we can apply the first main theorem of \cite{CM6}
 (theorem $0.1$ there) to get a subsequence $\tilde{\Sigma}_{i'}$
 converging off of a Lipschitz curve $\cS$ (where $|A| \to \infty$)
 to a foliation of $\RR^3$ by
 parallel planes.  Moreover, this Lipschitz curve goes through $0$
 and is transverse to the planes and consequently intersects every hemisphere above
 the plane through $0$.  However, this is a contradiction since  \eqr{e:c5}
   gives a   scale-invariant curvature bound  above
 this plane.
\end{proof}

Finally, we will need an easy consequence of the one-sided
curvature estimate of \cite{CM6} (this consequence is corollary
I.1.9 in \cite{CM6}):

\begin{Cor}   \label{c:conecor} \cite{CM6}
There exists   $\delta_0>0$ so that the following holds:

\noindent Let $\Sigma\subset B_{2R_0}$  be an embedded minimal
disk with $\partial \Sigma\subset
\partial B_{2R_0}$. If $\Sigma$  contains a $2$-valued graph $\Sigma_d \subset
\{x_3^2 \leq \delta_0^2\, (x_1^2+x_2^2)\}$ over $D_{R_0}\setminus
D_{r_0}$ with gradient $\leq \delta_0$, then each component of
\begin{equation}
    B_{R_0/2}\cap \Sigma \setminus (\cone_{\delta_0}(0)\cup B_{2
r_0}) \end{equation}
 is a multi-valued graph with gradient $\leq
1$.
\end{Cor}

\subsection{Properties (0)-(4) in Lemma \ref{l:t0.1}}

Properties (1) through (4) in Lemma \ref{l:t0.1} were implicit in
\cite{CM6} and we will describe below how to prove them using the
results in \cite{CM6}.

We first describe how to get the blow up points satisfying
(0)-(3).

\begin{itemize}
\item
{\bf{The slope $\delta$ and constant $C_1$}}: Set $\delta =
\delta_0$ from Corollary \ref{c:conecor}.  Then let $C_1$ and
$C_2$ be given by Lemma \ref{t:mv} with $N=2$ and $\epsilon =
\delta / 8$.
\item
{\bf{The initial multi-valued graph}}:
The lower curvature bound
\eqr{e:cmax} and a simple blowup argument (lemma $5.1$ in
\cite{CM4}) give a blow up pair $(y_0 , s_0)$ with
\begin{equation}
    B_{C_1 s_0}(y_0) \subset B_{C' \, R / c_{max}} \, .
\end{equation}
Lemma \ref{t:mv} then gives an associated rotation of $\RR^3$ and
a $2$-valued graph $\Sigma_0$ with gradient $\leq \delta / 8$ over
\begin{equation}
    D_{R/(2C_2)}(y_0) \setminus D_{C_1 \, s_0}(y_0) \, .
\end{equation}
  (Here we have used a slight abuse of notation since $y_0$ may
not be in the plane $\{ x_3 = 0 \}$.)
\item {\bf{Blow up pairs satisfying (0) are nearly parallel}}:
As long as $c_{out}$ is sufficiently large, then any blow up pair
$(y_i,s_i)$ satisfying (0) automatically has gradient $\leq \delta
/ 3$.  To see this, simply note that it has gradient $\leq \delta
/ 8$ over {\underline{some}} plane; embeddedness then forces this
plane to be almost parallel to the plane $\{ x_3 = 0 \}$.
\item {\bf{Nearby blow up pairs satisfy (0)}}:
After possibly choosing $c_{max}$ even larger,
then \eqr{e:cmax} implies that any blow up pair $(y_i , s_i)$ with
$y_i \in B_{2R/c_{out}}$ must have $C_1 \, s_i \leq 4R/c_{out}$,
i.e., must satisfy (0).
\item {\bf{Blow up pairs satisfying (0), (1), and (2)}}:  We will
iteratively apply Lemma \ref{l:nearby} to blow up pairs $(y_i ,
s_i)$ satisfying (0)--(2).  To get the first pair above $y_0$,
apply Lemma \ref{l:nearby} to get $(y_1 , s_1)$ above $y_0$ with
\begin{equation}    \label{e:bupin2}
    B_{C_1 \, s_{1}} (y_{1}) \subset B_{C_5 \, s_0}(y_0) \setminus
    B_{C_1 \, s_0}(y_0) \, .
\end{equation}
Repeat this to find $y_2$, etc., until
\begin{equation}
    B_{C_5 \, s_i} (y_i) \cap
    \partial B_{2R/c_{out}} \ne \emptyset \, .
\end{equation}
 The $y_i$'s with $i< 0$ are constructed similarly.  Note that
every $y_i$ is then contained in $B_{2R/c_{out}}$ so that (0)
holds.  Finally,  the cone property (1) follows immediately from
Corollary \ref{c:conecor}.
\item {\bf{Property (3)}}:  Iteratively applying (1) directly
gives (3).  This is because   (1) gives a lower bound for the
slope of the line segment connecting consecutive $y_i$'s.
\end{itemize}

We will next describe how to get (4) by combining (1)-(3) with
results of   \cite{CM3}-\cite{CM6}.  Finally, we will establish
(5) in the next subsection.

Observe first that Lemma \ref{t:mv}  directly gives the gradient
bound (4) on each of the corresponding $2$-valued graphs. To
extend this gradient bound to the rest of $\Sigma$, note that we
can choose a constant $C_2'$ so that each point
\begin{equation}
    y \in B_{R/C_2'} \cap \Sigma \setminus \cup_i B_{C_2' \, s_i}(y_i)
\end{equation}
satisfies a one-sided condition as in Corollary \ref{c:barrier}.
Precisely, $y$ is between  the $2$-valued graphs corresponding to
some $y_i$ and $y_{i+1}$ and, furthermore, these graphs are
themselves close enough together that we get two (in fact many)
distinct components  of
\begin{equation}
    B_{|y-y_i|/2} (y) \cap \Sigma
\end{equation}
which intersect the smaller concentric extrinsic ball
\begin{equation}
    B_{\epsilon \,|y-y_i|/(2c)}(y) \, .
\end{equation}
  Therefore, Corollary \ref{c:barrier} gives a curvature
estimate near $y$.  Finally, the desired gradient bound (4) at $y$
then follows from this curvature bound, the bound for the gradient
of the $2$-valued graphs $y$ is pinched between, and the gradient
estimate.  The fact that there are exactly two of these
multi-valued graphs was proven in proposition II.1.3 in
\cite{CM6}.

\subsection{The proof of (5) in Lemma \ref{l:t0.1}}

The key to establishing (5) is to first prove a chord arc bound
 assuming bounded
curvature (Lemma
\ref{l:ca}) and second to establish the curvature bound (Lemma
\ref{l:bound}).  This chord arc bound is essentially lemma
II.2.1 of \cite{CM6}, but the statement there does not suffice for
the application here. The statement that we need is the following:

\begin{Lem}  \label{l:ca}  (cf. lemma II.2.1 in \cite{CM6}.)
 There exists $C_d > 1$ so that given a constant $C_a$, we get another constant $C_b$
such that the following holds:

\noindent If $\Sigma \subset \RR^3$ is an embedded minimal disk
with $0 \in \Sigma \subset B_{R}$ and $\partial \Sigma \subset
\partial B_{R}$  and in addition
\begin{equation}        \label{e:cacb}
    \sup_{ B_{R} \cap \Sigma } |A|^2 \leq C_a \, R^{-2} \, ,
\end{equation}
 then
 \begin{equation}
    \Sigma_{0, \, \frac{R}{C_d} } \subset \cB_{C_b \, R}(0) \, .
\end{equation}
\end{Lem}

\begin{proof}
See Appendix \ref{s:a4}.
\end{proof}

The second result from \cite{CM6} that we will need is a curvature
bound on a
 larger extrinsic ball $B_{C_3 s_i}(y_i)$ around a blow up point $(y_i , s_i)$.
  The proof of this
curvature bound is essentially contained in the proof of lemma
I.1.10 in \cite{CM6} but was not made explicit there.  For
completeness, we state and prove this bound below:

\begin{Lem}    \label{l:bound}  \cite{CM6}
For every positive number $C_3$, there is a positive number $C_4$
with the following property.  If
\begin{enumerate}
\item[(a)]
$\Sigma \subset \RR^3$ is an embedded minimal disk with $\Sigma
\subset B_{C_4 \, s}(y)$ and $\partial \Sigma \subset
\partial B_{C_4 \, s}(y)$,
\item[(b)]
$(y,s)$ is a blow up pair,
\end{enumerate}
then we get the curvature bound
\begin{equation}   \label{e:c4}
    \sup_{B_{C_3 \, s}(y) \cap \Sigma} |A|^2 \leq C_4 \, s^{-2} \,
    .
\end{equation}
\end{Lem}

\begin{proof}
After rescaling and translating  $\Sigma$, we can assume
 that $y =0$ and $s =1$.
We will argue by contradiction; suppose therefore that $\Sigma_i$
 is a sequence of embedded minimal disks
  satisfying (a) and (b) with $y=0$, $s=1$, and $C_4 = i$
 but so that \eqr{e:c4} fails for some fixed $C_3$.

  Since both the
 radii $i$ of the extrinsic balls go to infinity  and
\begin{equation}
    \sup_{B_{C_3}(0) \cap \Sigma_i} |A|^2 \to \infty \, ,
\end{equation}
  we can apply the first main theorem of \cite{CM6}
 (theorem $0.1$ there).  Therefore, a subsequence $\Sigma_{i'}$
 converges off of a Lipschitz curve $\cS$ to a foliation of $\RR^3$ by
 parallel planes.  This convergence implies that the supremum of
 $|A|^2$ on each fixed extrinsic ball either goes to zero or
 infinity, depending on whether or not this ball intersects $\cS$.
 However, this directly
 contradicts the assumption (b), thereby giving the lemma.
\end{proof}

To prove (5), we first use Lemma  \ref{l:bound} to get a uniform
curvature bound on larger extrinsic balls $B_{C_3' s_i}(y_i)$.
Combining Lemma \ref{l:ca},  and using the one-sided estimate
(i.e., Corollary \ref{c:barrier}) to see that there is only such
component,  then gives (5).

\subsection{The proof of Proposition \ref{t:c-a}}   \label{ss:25}

We will next see how properties (0)--(5) in Lemma \ref{l:t0.1}
imply Proposition \ref{t:c-a}.

\begin{proof}
(of Proposition \ref{t:c-a}).  We will divide the proof into two
cases, depending on whether or not the curvature is large, i.e.,
whether \eqr{e:cmax} holds.

Suppose first that \eqr{e:cmax} fails so that   we have the
curvature bound
\begin{equation}    \label{e:notcmax}
    \sup_{B_{  R /c_{max} }(x) \cap \Sigma} |A|^2 \leq c_{max}^2 \, R^{-2} \,
    .
\end{equation}
We can then  apply Lemma \ref{l:ca} to get
\begin{equation}
    \Sigma_{x,c_1' \, R} \subset \cB_{ c_1 R }(x) \, ,
\end{equation}
giving the proposition in this case.

 In the second case, where \eqr{e:cmax} holds, the
proposition will follow from Lemma \ref{l:t0.1}.  We do this in
two steps.

First, for any point
\begin{equation}
    z \in B_{\delta R/
    (4c_{out})} (x) \cap \Sigma \, ,
\end{equation}
 we have
\begin{equation}    \label{e:leafd}
    \dist_{\Sigma} \, ( z , \,  \cup_{i} B_{c_{in} \, s_i}(y_i) ) \leq
    C' \,
    R  \, .
\end{equation}
 This follows immediately from the
gradient bound for the multi-valued graphs given by (4) together
with the fact that the points $y_i$ go from top to bottom by (2)
and (3).

Second, (1) and (5) imply a bound for the diameter of the union of
the balls $B_{c_{in} \, s_i}(y_i)$.  Namely,   the balls
$B_{s_i}(y_i)$ are disjoint and satisfy the cone property (1) and,
therefore, we get a bound for the sum of the radii $s_i$ of these
balls
\begin{equation}
    \sum_i s_i \leq C_0 \, R/c_{in} \, .
\end{equation}
 Combining this with the chord arc property (5) then
gives a bound for the diameter of the union of these balls
\begin{equation}    \label{e:axisd}
    \diam_{\Sigma} \, ( B_{R/c_{out}}(x) \cap \cup_{i} B_{c_{in} \, s_i}(y_i) )   \leq C' \, R  \, .
\end{equation}
  Combining the   bounds \eqr{e:leafd} and \eqr{e:axisd}, the triangle
inequality
  gives the proposition in this case as well.
\end{proof}

\section{The proof of Proposition \ref{t:2}}   \label{s:s3}

 In this section, we will complete the proof of  Proposition \ref{t:2}.
 To do this, we will first define a   weak chord arc property for an
intrinsic ball.  This  property requires that the intrinsic ball
contains an entire component of $\Sigma$ in a smaller extrinsic
ball.

 Throughout this section  $\Sigma \subset \RR^3$ is an
embedded minimal disk, possibly noncompact, with boundary
$\partial \Sigma$.

\subsection{Weakly chord arc}

To show Proposition \ref{t:2}, we need to prove that there is a
constant $\delta_1
> 0$ so that for any intrinsic ball $\cB_R (x) \subset \Sigma \setminus
\partial \Sigma$ we have the inclusion
\begin{equation}    \label{e:deflp}
    \Sigma_{x,\delta_1 \, R} \subset \cB_{ \frac{R}{2} }(x) \, ,
\end{equation}
where, as before, $\Sigma_{x,\delta_1 \, R}$ denotes the component
of $B_{\delta_1 \, R}(x) \cap \Sigma$ containing $x$.

 Since $\Sigma$ is smooth,  the inclusion
\eqr{e:deflp} must hold for  sufficiently small balls depending on
$\Sigma$.  The key step in the proof of Proposition \ref{t:2} is
to show that if \eqr{e:deflp} holds on one scale, then it also
holds on five times the scale.  (Here, when we say that it holds
on a scale, we mean that it holds for all balls of this radius;
cf. (A') in the proof.)  This will be done in Proposition
\ref{p:key} below.   Proposition \ref{t:2} will then follow by
using a  blow up argument (Lemma \ref{l:scaling} below) to locate
the largest scale where \eqr{e:deflp} holds and then applying
 Proposition \ref{p:key} to see that  \eqr{e:deflp}
  continues to hold on larger scales.

 We will say that an intrinsic ball where we have the inclusion  \eqr{e:deflp}
 is weakly chord arc; namely, we make the following definition:

\begin{Def} (Weakly chord arc).
An intrinsic ball $\cB_s(x) \subset \Sigma \setminus \partial
\Sigma$ is said to be {\it{$\delta$-weakly chord arc}} for some
$\delta
> 0$  if \eqr{e:deflp} holds with $R=s$ and $\delta = \delta_1$.
  Note that
\eqr{e:deflp} is only possible if $\delta \leq 1/2$.
\end{Def}

It will be important later that subballs of a weakly chord arc
ball are themselves weakly chord arc.  While this does not follow
directly from \eqr{e:deflp}, we do directly get that the
intersections with smaller extrinsic balls are compact and have
boundary in the boundary of the smaller ball.  In particular,
these properties will allow us to apply Proposition \ref{t:c-a} to
conclude that the smaller balls are themselves $\delta_2$-weakly
chord arc; this will be done in the beginning of the proof of
Proposition \ref{t:2} when we replace (A) with (A') there.

It will be convenient to introduce notation for the largest radius
of a weakly chord arc ball about a given point.  We will do this
next.

Given a constant $\delta$ and a point $x \in \Sigma \setminus
\partial \Sigma$, we let
 $\rdel (x)$ denote the largest radius where $\cB_{\rdel (x)}(x)$ is
 $\delta$-weakly chord arc, i.e.,
 \begin{equation}   \label{e:defrx}
    \rdel (x) = \sup \, \{ R \, | \, \cB_{R}(x) \subset
    \Sigma \setminus \partial \Sigma {\text{ is $\delta$-weakly
    chord arc }} \} \, .
 \end{equation}
  Since $\Sigma$ is a smooth surface,
 we obviously have $\rdel (x) > 0$ for every $x$ and any $\delta < 1/2$.

We can now state the key proposition which shows that if all
intrinsic balls of radius $R_0$ near a point $y$ are
$\delta_2$-weakly chord arc, then so is the five-times ball
$\cB_{5\,R_0}(y)$  about $y$. The constant $\delta_2$ in the
proposition is given by Proposition \ref{t:c-a}.

\begin{Pro}     \label{p:key}
Let $\Sigma \subset \RR^3$ be an embedded minimal disk. There
exists a constant $C_b
> 1$ independent of $\Sigma$ so that if
 $\cB_{C_b \,
R_0}(y) \subset \Sigma \setminus \partial \Sigma$ is an intrinsic
ball and
\begin{enumerate}
  \item[(A')]
    every intrinsic subball $\cB_{R_0}(z) \subset \cB_{ C_b \,R_0}(y)$
    is $\delta_2$-weakly chord arc,
\end{enumerate}
then, for every $s \leq 5\, R_0$, the intrinsic ball $\cB_{s}(y)$
    is $\delta_2$-weakly chord arc.
\end{Pro}

\subsection{Extrinsically close yet intrinsically far apart}

In this subsection, we recall from \cite{CM2} and \cite{CM4}
several important properties of embedded minimal surfaces with
bounded curvature. The basic point is that nearby, but disjoint,
minimal surfaces with bounded curvature can be written as graphs
over each other of a positive function $u$ which satisfies a
useful second order elliptic equation.  We will focus here on two
consequences of this.  The first is a chord arc result assuming an
a priori curvature bound (see Lemma \ref{l:halfstable} below). The
second is that this elliptic equation for $u$ implies a Harnack
inequality for $u$ that bounds the rate at which the two disjoint
surfaces can pull apart.

 We will need the
notion of $1/2$-stability. Recall from \cite{CM4} that a domain
$\Omega \subset \Sigma$ is said to be $1/2$-stable
 if, for
all Lipschitz functions $\phi$ with compact support in $\Omega$,
  we have the $1/2$-stability
inequality:
\begin{equation} \label{e:stabineq}
    \frac{1}{2} \, \int |A|^2 \, \phi^2
        \leq  \int |\nabla \phi |^2 \,  .
\end{equation}

Loosely speaking, the next elementary lemma shows that if two
disjoint intrinsic balls are extrinsically close (see
\eqr{e:delta}) and have a priori curvature bounds (see
\eqr{e:ca}), then smaller concentric intrinsic balls are almost
flat and thus in particular their boundaries are far away from
their centers (see the conclusion \eqr{e:conclu}). Since it is
only this last conclusion that we need, and not the stronger
statement that they are almost flat, we only state this.

\begin{Lem}     \label{l:halfstable}
There exists $C_0 > 1$ so that for every $C_a > 0$, there exists
  $\tau > 0$ such that if $\cB_{C_0} (x_1)$ and $\cB_{C_0}
(x_2)$ are disjoint intrinsic balls in $\Sigma \setminus \partial
\Sigma$ with
\begin{align}   \label{e:ca}
    \sup_{\cB_{C_0} (x_1) \cup \cB_{C_0} (x_2)} |A|^2 &\leq C_a \,
    , \\
    |z_1 - z_2| &< \tau \, , \label{e:delta}
\end{align}
then for $i=1$, $2$ we have
\begin{equation}    \label{e:conclu}
    B_{10} (x_i) \cap \partial \cB_{11} (x_i) = \emptyset \, .
\end{equation}
\end{Lem}

\begin{proof}
 Using the argument of
\cite{CM2} (i.e., curvature estimates for $1/2$-stable surfaces)
we get a constant
 $C_0 > 1$ so that if
$\cB_{C_0/2 }(z) \subset  \Sigma \setminus \partial \Sigma$ is
$1/2$-stable, then $\cB_{11}(z)$ is a graph with
\begin{equation}
    B_{10}(z) \cap \partial \cB_{11}(z) = \emptyset \,
 .
 \end{equation}
 Corollary 2.13 in \cite{CM4} gives $\tau = \tau (C_a)
> 0$ so that if $ |z_{1} - z_{2}| < \tau$ and $|A|^2 \leq C_a$
on (the disjoint balls) $\cB_{ C_0 }(z_i)$, then each subball
\begin{equation}
    \cB_{\frac{C_0 }{2} } (z_i) \subset \Sigma
\end{equation}
     is $1/2$-stable.
\end{proof}

As mentioned above, one of the key points in the proof of the
previous lemma was that nearby, but disjoint, embedded minimal
surfaces with bounded curvature can be written as graphs over each
other of a positive function $u$. Furthermore, standard
calculations show that this function $u$ satisfies a second order
elliptic equation resembling the Jacobi equation (for the Jacobi
equation, the functions $a_{ij} , b_j , c$  in \eqr{e:wlog2}
vanish). These standard, but very useful, calculations were
summarized in lemma $2.4$ of \cite{CM4} which we recall next.

\begin{Lem} \label{l:mge} \cite{CM4}
There exists $\delta_g > 0$ so that if $\Sigma$ is minimal and if
$u $ is a positive solution of the minimal graph equation over
$\Sigma$ (i.e., $\{ x + u(x) \, \nn_{\Sigma} (x) \, |\, x \in
\Sigma \}$ is minimal) with
\begin{equation}    \label{e:mge}
    |\nabla u| + |u| \,
|A| \leq \delta_g \, ,
\end{equation}
 then $u$
satisfies on $\Sigma$
\begin{equation} \label{e:wlog2}
\Delta u =  \dv (a \nabla u)  + \langle b , \nabla u \rangle +
(c-1) |A|^2 \, u \, ,
\end{equation}
for functions $a_{ij} , b_j , c$ on $\Sigma$ with $|a| , |c| \leq
3 \, |A| \, |u| + |\nabla u|$ and $|b| \leq 2 \, |A| \, |\nabla
u|$.
\end{Lem}

Equation \eqr{e:wlog2} implies a uniform Harnack inequality for
$u$ which bounds the supremum of $u$ on a compact subset of
$\Sigma \setminus \partial \Sigma$ by a multiple of the infimum;
see, for instance, theorem $8.20$ in \cite{GiTr}. We will use this
in the next subsection to show that two nearby, but disjoint,
components of $\Sigma$ with bounded curvature pull apart very
slowly.

\subsection{Extending weakly chord arc to a larger scale: The proof of Proposition \ref{p:key}}

We are now prepared to prove Proposition \ref{p:key}, i.e., to
show that if all intrinsic balls of radius $R_0$ near a point $y$
are weakly chord arc, then so is the five-times ball
$\cB_{5\,R_0}(y)$  about $y$.
 To do this, we first show that $\cB_{5\,R_0}(y)$
 is still weakly chord arc, but with a
worse constant. We then use Proposition \ref{t:c-a} to improve the
constant, i.e., to see that it is in fact $\delta_2$-weakly chord
arc.

The reader may find it helpful to compare the proof below with the
  simpler proof of the special case where $\Sigma$ has bounded
curvature, i.e., with the proof of Lemma \ref{l:ca} given in
Appendix \ref{s:a4}.  The difference is that here the one-sided
curvature estimate is used, while there we simply assume an a
priori bound on the curvature.

\begin{proof}
(of Proposition \ref{p:key}.)  After rescaling  and translating
$\Sigma$, we can assume that $R_0 = 1$ and $y=0$.

 The proposition follows from the following claim:
There exists $n$ so that
 \begin{equation}        \label{e:claim}
        \Sigma_{0,\, 5} \subset
        \cB_{(6n+3) \, C_0}(0) \, ,
\end{equation}
where $C_0 > 1$ is given by Lemma \ref{l:halfstable}.  The
proposition will follow immediately from \eqr{e:claim} by applying
Proposition \ref{t:c-a} to $\Sigma_{0,\, 5}$.  Namely,
\eqr{e:claim} implies that the embedded minimal disk $\Sigma_{0,\,
5}$ is compact and has
\begin{equation}
    \partial \Sigma_{0,\, 5} \subset  \partial B_5 \, .
\end{equation}
We can therefore apply Proposition \ref{t:c-a} for any $t \leq 5$
to get that
\begin{equation}   \label{e:c-a2}
    \Sigma_{0, \, \delta_2 \, t} \subset \cB_{t/2}(0) \, ,
\end{equation}
giving the proposition.

\vskip2mm We will prove
   the claim (i.e., \eqr{e:claim}) by arguing by contradiction; so suppose that
 \eqr{e:claim}
fails for some
large $n$.  Consequently, we get a curve
\begin{equation}
    \sigma \subset \Sigma_{0,\, 5} \subset B_{5}
\end{equation}
 from $0$ to a point in $\partial \cB_{(6n+3)\, C_0}(0)$. For $i=1, \dots ,
n$, fix points
\begin{equation}
    z_i \in \partial \cB_{6i\,  C_0}(0) \cap \sigma \, .
\end{equation}
  It follows that the intrinsic balls $\cB_{3\, C_0}(z_i)$:
\begin{itemize} \item
 Are
disjoint. \item Have centers in $B_{5} \subset \RR^3$.
\end{itemize}

Since the $n$ points $\{ z_i \}$ are all in the Euclidean ball
$B_5 \subset \RR^3$,   there exist integers $i_1$ and $ i_2$ with
\begin{equation}    \label{e:427}
    0< | z_{i_1} - z_{i_2}| < C' \,  n^{-1/3}  \, .
\end{equation}
   Furthermore, since each intrinsic ball of
radius one about any $z_i$ is $\delta$-weakly chord arc by (A'),
we have that each embedded minimal disk $\Sigma_{z_i, \delta}$ is
compact and has
\begin{equation}     \label{e:428}
    \partial \Sigma_{z_i, \delta} \subset \partial B_{\delta}(z_i)
    \, .
\end{equation}
 Consequently, for $n$ large enough, \eqr{e:427} implies that
 the components $\Sigma_1$ and $\Sigma_2$ of
 \begin{equation}
    B_{\frac{\delta}{2}}(z_{i_1}) \cap \Sigma
 \end{equation}
 containing $z_{i_1}$ and $z_{i_2}$, respectively, are compact and
 have
\begin{equation}     \label{e:429}
    \partial \Sigma_{i} \subset \partial B_{\frac{\delta}{2}}(z_{i_1})
    \, .
\end{equation}
Note that the center of this extrinsic ball is the same for
$\Sigma_1$ and $\Sigma_2$.  Let $c>1$ be given by Corollary
\ref{c:barrier}.
 For $n$ sufficiently large,
   \eqr{e:427} implies that   $\Sigma_2$
   intersects the smaller concentric extrinsic ball
   $B_{\frac{\delta}{2c}}(z_{i_1})$ and, since
     $\Sigma_1$ contains the center of this ball, then it follows
   that for both $j=1$ and $j=2$, we have that
\begin{equation}    \label{e:430}
    B_{\frac{\delta}{2c}}(z_{i_1}) \cap \Sigma_j \ne \emptyset \,
    .
\end{equation}
   Combining
\eqr{e:429} and \eqr{e:430}, Corollary \ref{c:barrier} gives the
curvature bound
  for $j=1$, $2$
\begin{equation}    \label{e:cafromb}
    \sup_{ \cB_{ \frac{\delta}{2c} }(z_{i_j}) } \, |A|^2 \leq \left(
    \frac{\delta}{2c}\right)^{-2}
    \, .
\end{equation}

By lemma $2.11$ of \cite{CM4}, the curvature bound \eqr{e:cafromb}
gives a constant $r' = r' (\delta , c)$ so that if $n$ is
sufficiently large, then
    $\cB_{ 3 \, r' }(z_{i_2})$  can be written as a normal exponential graph of a function $u$
    over a domain
    $\Omega$, where:
\begin{enumerate}
\item[(i)]
The function $u$ satisfies \eqr{e:mge}. \item[(ii)] The domain
$\Omega$ contains, and is contained in, concentric intrinsic balls
as follows
\begin{equation}
    \cB_{ 2 \, r' }(z_{i_1}) \subset \Omega
        \subset \cB_{ 4 \, r' }(z_{i_1})  \, .
\end{equation}
\end{enumerate}
(To see this, first use the curvature bound to write each
component locally as a graph and then use embeddedness to see that
these graphs must be roughly parallel.)
   By Lemma \ref{l:mge} (and \eqr{e:427}), we
can apply the Harnack inequality to $u$ to get
\begin{equation}    \label{e:indu1}
    \sup_{  \cB_{ r'  }(z_{i_1}) }
         u \leq \tilde{C} \, |z_{i_2} -
    z_{i_1}| \leq \tilde{C}' \, n^{-1/3} \, .
\end{equation}
As long as $n$ is large enough, \eqr{e:indu1} allows us to repeat
the argument with a point in the boundary $\partial
\cB_{r'}(z_{i_1})$ in place of  $z_{i_1}$.  Therefore, for
 $n$   large enough, we can repeatedly combine Corollary \ref{c:barrier} and
  the Harnack
 inequality to extend the curvature bound \eqr{e:cafromb} to the larger intrinsic
 balls
 \begin{equation}   \label{e:theballs}
    \cB_{C_0}(z_{i_j}) {\text{ for }} j=1, \, 2 \, .
 \end{equation}

 Now that we have a uniform curvature bound on the disjoint intrinsic balls
 \eqr{e:theballs} and the centers of these balls are extrinsically close by
 \eqr{e:427}, we can apply Lemma \ref{l:halfstable} to get that
\begin{equation}
    B_{5} \cap
\partial \cB_{11}(z_{i_j})  = \emptyset \, .
\end{equation}
(Here we used that $B_{5}  \subset B_{10}(z_{i_j})$ because
$z_{i_j} \in B_5$.)
 Since the curve $\sigma$
must intersect $\partial \cB_{11}(z_{i_j})$, this contradicts
 the fact that the curve $\sigma$ is contained in
the ball $B_{5}$.  This  contradiction proves \eqr{e:claim} and
  gives the proposition.
\end{proof}

The previous proposition is the key step in the proof of
Proposition \ref{t:2}.  To complete the proof, we will use a
simple blow up argument to find some small initial scale which is
weakly chord arc and then apply Proposition \ref{t:2} to get that
so are larger scales.   As is often the case in this type of blow
up argument, then the existence of such an initial scale is
complicated slightly by the fact that $\Sigma$ has non-empty
boundary.

To incorporate the   boundary, we let
  $a_{\delta}$  be the supremum of the ratio of the
distance to $\partial \Sigma$ to the largest radius of an
intrinsic ball which is $\delta$-weakly chord arc, i.e., set
\begin{equation} \label{e:s2}
    a_{\delta} = \sup_{ z \in \Sigma } \, \,  \frac{\dist_{\Sigma} (z , \partial \Sigma)}{
        \rdel (z)  }\,
        ,
\end{equation}
where $\rdel (z)$ is given by \eqr{e:defrx}.

\subsection{Upper bounds for $a_{\delta}$}

Suppose for a moment that $\Sigma$ is compact and smooth up to the
boundary $\partial \Sigma$ and $\delta < 1/2$.  We will, in the
proof of Lemma \ref{l:scaling} below, use  that
\begin{equation}    \label{e:adeltaf}
    a_{\delta} < \infty \, .
\end{equation}
  To see \eqr{e:adeltaf},
observe that compactness and smoothness give uniform bounds on
$|A|^2$ and the geodesic curvature of $\partial \Sigma$.  Given
any constant $\epsilon > 0$, the bound on $|A|^2$ gives a constant
$r_0 > 0$ so that if $s \leq r_0$ and $\cB_{s}(z) \subset \Sigma
\setminus \partial \Sigma$, then $\cB_s (z)$ is a graph over some
plane of a function with gradient $\leq \epsilon$. In particular,
the intrinsic ball $\cB_s (z)$ is $\delta$-weakly chord arc for
$\epsilon$ sufficiently small. Furthermore, the bound on the
geodesic curvature of $\partial \Sigma$ gives a constant $r_1 > 0$
so that if
\begin{equation}    \label{e:dbound}
    d_z = \dist_{\Sigma}(z, \partial \Sigma) \leq r_1 \, ,
\end{equation}
 then $\cB_{d_z}(z) \subset \Sigma \setminus \partial \Sigma$.
We can then establish \eqr{e:adeltaf} by considering two cases
depending on the distance to the boundary.  If
\begin{equation}    \label{e:dbound2}
    d_z = \dist_{\Sigma}(z, \partial \Sigma) \leq \min \, \{ r_0 , r_1 \} \, ,
\end{equation}
 then $\cB_{d_z}(z)$ is $\delta$-weakly chord arc so that
 \begin{equation}
    \rdel (z)
 = \dist_{\Sigma} (z , \partial \Sigma) \, .
 \end{equation}
   On the other hand, when
 \eqr{e:dbound2} fails, then $\cB_{r_2}(z)$ is $\delta$-weakly chord arc
 where $r_2 = \min \, \{ r_0 , r_1 \}$ and hence
 \begin{equation}
   \frac{\dist_{\Sigma} (z , \partial \Sigma)}{
        \rdel (z)  } \leq  \frac{\diam \, (\Sigma) }{
        r_2  } \, .
 \end{equation}
This shows that $a_{\delta} < \infty$ if $\Sigma$ is compact and
smooth.

Let us return to Proposition \ref{t:2}.   It is not hard to see
that the proposition   is equivalent to
 an upper bound (independent of $\Sigma$)
  for $a_{\delta}$ for a fixed $\delta > 0$.
 Namely, suppose that $\cB_R(x) \subset \Sigma \setminus \partial \Sigma$
 is as in the proposition and
we have an upper bound for $a_{\delta}$
\begin{equation}    \label{e:cdelta}
    a_{\delta} \leq c  < \infty \, .
\end{equation}
  Since $\cB_R(x) \subset \Sigma
\setminus \partial \Sigma$, then \eqr{e:s2}
 implies that
\begin{equation}    \label{e:defit}
     {R}  \leq {\dist_{\Sigma} (x , \partial \Sigma)} \leq c \,
     {
        \rdel (x)  }  \, .
\end{equation}
   Consequently, by the definition
\eqr{e:defrx} of $\rdel (x)$,   there exists a radius $s
> \frac{R}{2\, c}$ so that $\cB_s(x)$ is $\delta$-weakly
chord arc and hence
\begin{equation}    \label{e:tgives}
    \Sigma_{ x , \, \frac{\delta\, R}{4 \, c } } \subset
    \Sigma_{x, \frac{\delta \, s}{2} } \subset
    \cB_{\frac{R}{2}}(x) \, .
\end{equation}
Equation   \eqr{e:tgives} would
 then give Proposition \ref{t:2}.

\subsection{Locating the smallest scale which is not weakly chord arc}

 We will first need to locate a smallest scale on which $\Sigma$ is
 not $\delta$-weakly chord arc.  We do this in the next lemma with a
 simple blow up argument.  The $\Sigma$
  in this lemma is assumed to be compact and smooth up to the boundary
   so that $a_{\delta} <
  \infty$ by \eqr{e:adeltaf}.

\begin{Lem}     \label{l:scaling}
  Given $\Sigma$ compact and smooth up to the boundary
  and a constant $\delta$ with $0 < \delta < 1/2$,
  there exists $y \in \Sigma$ and $R_0>0$ so that:
  \begin{enumerate}
  \item[(A)]
     $\rdel (x) > R_0$ for every $x \in \cB_{ a_{\delta}
    \,R_0}(y)$,
    where $\rdel (x)$ is  given by \eqr{e:defrx}.
\item[(B)]
    The intrinsic ball $\cB_{5 \, R_0}(y)$ is {\underline{not}} $\delta$-weakly chord arc.
  \end{enumerate}
\end{Lem}

\begin{proof}
  Define a function $G$ on $\Sigma$ by setting
\begin{equation}
    G(x) = \frac{ \dist_{\Sigma} (x , \partial
    \Sigma)}{\rdel (x)} \, .
\end{equation}
  Since $\Sigma$ is smooth and compact, then
\eqr{e:adeltaf} and the definitions of $G$ and $a_{\delta}$  give
that
\begin{equation}
      a_{\delta} = \sup G < \infty \, .
\end{equation}
  We
can therefore choose $y$ so that $G(y)$ is greater than half the
supremum $a_{\delta}$ of $G$ on $\Sigma$.  Hence,
\begin{equation}    \label{e:s3}
    \frac{ \dist_{\Sigma} (y,\partial \Sigma) }{\rdel (y)} =
    G(y ) >   \frac{\sup G }{ 2} = \frac{ a_{\delta} }{ 2} \, .
\end{equation}
We will see that \eqr{e:s3} implies (A) and (B) with $R_0 = \rdel
(y) / 4$.

Set $d_{\partial} = \dist_{\Sigma} (y, \partial \Sigma)$ so that
if $x \in \cB_{d_{\partial}  /2}(y)$, then by the triangle
inequality
\begin{equation}    \label{e:s4}
     \dist_{\Sigma} (x , \partial
    \Sigma) >  \frac{ d_{\partial} }{ 2 }  \, .
\end{equation}
Combining \eqr{e:s3} and \eqr{e:s4} gives for $x \in
\cB_{d_{\partial} /2}(y)$ that
\begin{equation}
    \frac{ d_{\partial} }{ 2 \, \rdel (x) } < G(x) < 2 \, G(y)
        =  \frac{ 2 \, d_{\partial} }{ \rdel (y) } \, ,
\end{equation}
and thus
\begin{equation}    \label{e:co1}
     \rdel (x) > \frac{\rdel (y)}{4} = R_0 \, .
\end{equation}
  From \eqr{e:s3}, we see that   $2 \, a_{\delta} \, R_0 <
d_{\partial}$ and hence
\begin{equation}    \label{e:co2}
    \cB_{ a_{\delta} {R_0} }
    (y) \subset \cB_{ \frac{d_{\partial} }{2} }(y) \, .
\end{equation}
Combining \eqr{e:co1}  and  \eqr{e:co2} gives (A).  We get (B)
immediately from the maximality of $\rdel (y)$.
\end{proof}

\subsection{The proof of Proposition \ref{t:2}: Bounding $a_{\delta}$}

We are now prepared to prove Proposition \ref{t:2}, i.e., to show
that sufficiently small intrinsic balls in $\Sigma$ are weakly
chord arc.
 As mentioned above, this is equivalent to giving a
uniform upper bound for the constant $a_{\delta}$ defined in
\eqr{e:s2} for some fixed $\delta > 0$ (the constant $\delta$ will
be given by Proposition \ref{t:c-a}). In the actual proof, we will
first use Lemma \ref{l:scaling} to find the smallest scale which
is not $\delta$-weakly chord arc. To bound $a_{\delta}$, it
suffices to give a lower bound for this scale in terms of the
distance to the boundary $\partial \Sigma$. This is precisely the
content of Proposition \ref{p:key}.

\begin{proof}
(of Proposition \ref{t:2}.)
 Let the constant $\delta=\delta_2$ be
given by Proposition \ref{t:c-a}.   As we have seen in
\eqr{e:tgives},   the proposition follows from a uniform upper
bound for the constant $a_{\delta}$ defined in \eqr{e:s2}. The
rest of the proof is to establish such a bound.

Apply first Lemma \ref{l:scaling} to   locate the smallest scale
which is not $\delta$-weakly chord arc.  This gives a point $y$ in
$\Sigma$ and an intrinsic ball
 $\cB_{a_{\delta} \,
R_0}(y)$ so that:
\begin{enumerate}
  \item[(A)]
  $\rdel (z) > R_0$ for every $z \in \cB_{ a_{\delta}
    \,R_0}(y)$.
  \item[(B)]
  $\cB_{5{R_0}}(y)$ is {\underline{not}} $\delta$-weakly chord arc.
  \end{enumerate}
The condition (A) implies that each point $z \in \cB_{ a_{\delta}
    \,R_0}(y)$ is the center of
some $\delta$-weakly chord arc intrinsic ball of radius greater
than $R_0$.  However,  Proposition \ref{t:c-a} then easily gives
that $\cB_{R_0} (z)$ is in fact $\delta$-weakly chord arc (here we
use that $\delta$ is given by that proposition).  Namely, (A) can
be replaced by:
\begin{enumerate}
  \item[(A')]
    Every intrinsic ball $\cB_{R_0}(z)$ with $z \in \cB_{ a_{\delta} \,R_0}(y)$
    is $\delta$-weakly chord arc.
\end{enumerate}
The proposition now follows from Proposition \ref{p:key}.  Namely,
Proposition \ref{p:key} gives a constant $C_b$ so that if
\begin{equation}    \label{e:adelbig}
    a_{\delta} \geq C_b \, ,
\end{equation}
then (A') implies that the five times intrinsic ball
$\cB_{5{R_0}}(y)$ is $\delta$-weakly chord arc.  Since this would
contradict (B), we conclude that \eqr{e:adelbig} cannot hold and
the proposition follows.
\end{proof}

\section{Finite topology: The proofs of Corollaries \ref{c:finite}
and \ref{c:3rd}}

In this section, we prove both of Calabi's conjectures and
properness for complete embedded minimal surfaces with finite
topology. Recall that a surface $\Sigma$ is said to have finite
topology if it is homeomorphic to a closed Riemann surface of
genus $g$ with a finite set of  punctures. Each puncture
corresponds to an end of $\Sigma$ and thus the ends can be
represented by   punctured disks, i.e., each end is homeomorphic
to the set
\begin{equation}
    \{ z \in \CC \, | \, 0 < |z| \leq 1 \} \, .
\end{equation}

\subsection{Simply connected outside a compact set}

The key point for extending our results to surfaces with finite topology is to
show that intrinsic balls are eventually simply connected so that
our results for disks can be applied.  This is made precise in the
next lemma.

\begin{Lem} \label{l:ftc}
Let $\Gamma$ be a complete noncompact embedded minimal annulus
which contains one compact component $\gamma$ of $\partial
\Gamma$;  the other boundary is at infinity.  There is a constant
$\bar{R}$ (depending on $\Gamma$) so that the following holds:
\begin{equation} \label{e:(D)}
    {\text{If  $d_x = \dist_{\Gamma} (x , \gamma)  > \bar{R}$, then
    the intrinsic ball $\cB_{d_x/2}(x)$ is a disk.}}
\end{equation}
\end{Lem}

\begin{proof}
Suppose that \eqr{e:(D)} fails for every $\bar{R}$. It will follow
from that   $\Gamma$ is an annulus with non-positive curvature
that $\Gamma$ has finite total curvature. Namely, if \eqr{e:(D)}
fails, we get
 a sequence $x_i \in \Gamma$ with
 \begin{equation}  \label{e:di}
    d_i = \dist_{\Gamma} (x_i , \gamma) \to \infty
\end{equation}
 so that
the exponential map from $x_i$ is not injective into
$\cB_{d_i/2}(x_i)$.  In particular, there are distinct geodesics
$\gamma_i^a$ and $\gamma_i^b$ in $\cB_{d_i/2}(x_i)$ from $x_i$ to
a point $y_i \in \cB_{d_i/2}(x_i)$ and the closed curve
\begin{equation}
    \gamma_i = \gamma_i^a \cup \gamma_i^b
\end{equation}
is homologous to the compact boundary component $\gamma$.  Let
$\Gamma_i$ be the bounded component of $\Gamma\setminus \gamma_i$;
so $\Gamma_i$ is topologically an annulus bounded by $\gamma$ and
the piecewise smooth closed geodesic $\gamma_i$ with breaks at
$x_i$ and $y_i$.  Write $\int_{\gamma}k_g$ and
$\int_{\gamma_i}k_g$ for the two boundary terms in the
Gauss-Bonnet theorem for the annulus $\Gamma_i$ (both are
uniformly bounded; $\int_{\gamma_i}k_g$ is afterall just the angle
contribution at $x_i$ and $y_i$).  It follows that
\begin{equation}    \label{e:gi}
    \int_{\Gamma_i} |A|^2 = -2 \, \int_{\Gamma_i} K_{\Gamma} =  2
    \int_{\gamma} k_g + 2 \int_{\gamma_i} k_g \leq C \, .
\end{equation}
Moreover, by the triangle inequality, we have that $\dist_{\Gamma}
(\gamma , \gamma_i) \geq d_i / 2$ and hence $\Gamma_i$ contains
the intrinsic $(d_i/2)$-tubular neighborhood of $\gamma$.  Since
$d_i \to \infty$, the $\Gamma_i$'s exhaust $\Gamma$, i.e., $\Gamma
\subset \cup_i \Gamma_i$, and thus \eqr{e:gi} implies that
$\Gamma$ has finite total curvature.

Finally, we will show that \eqr{e:(D)} must hold when $\Gamma$ has
finite total curvature.  To see this, note that since $\Gamma$ is
an embedded annulus with finite total curvature, it is asymptotic
to either a plane or half of a catenoid (see, e.g., \cite{Sc2}).
In either case, \eqr{e:(D)} must hold for points sufficiently far
from the interior boundary $\gamma$.  This completes the proof of
the lemma.
\end{proof}

\subsection{Compact embedded annuli in a halfspace}

We will next bound the total curvature for a compact embedded
minimal annulus in a halfspace.  In the next lemma, we will use
$\Gamma_{\gamma , R}$ to denote the component  of $B_R \cap
\Gamma$ containing the boundary component $\gamma$.

\begin{Lem}     \label{l:annuli1}
Let $\Gamma$ be as in Lemma \ref{l:ftc}.  There exist constants
$\epsilon > 0$ and $\hat{R}$ so that if $R
> \hat{R}$, the component $\Gamma_{\gamma , 2R}$ is
compact, and
\begin{equation}
    \Gamma_{\gamma , 2R} \subset  \{ x_3 > - \epsilon \, R \} \, ,
\end{equation}
then  $\Gamma_{\gamma , R}$   has bounded total curvature
\begin{equation}        \label{e:btc}
    \int_{\Gamma_{\gamma , R}} |A|^2 \leq 2 \, \int_{\gamma}
    k_g + 8\pi \, .
\end{equation}
\end{Lem}

\begin{proof}
The bound \eqr{e:btc} follows immediately from the Gauss-Bonnet
theorem and the following two claims:
\begin{enumerate}
\item[(C1)]
There is a constant $\epsilon > 0$ so that if $\Gamma_{\gamma ,
2R} \subset \{ x_3 > - \epsilon \, R \}$ and $\partial
\Gamma_{\gamma , R} \setminus \gamma$ intersects $\{ x_3 <
\epsilon \, R \}$, then $\partial \Gamma_{\gamma , R} \setminus
\gamma$ is a graph over (a curve in) $\{ x_3 = 0 \}$ and
\begin{equation}    \label{e:cbd}
    \int_{\partial \Gamma_{\gamma , R} \setminus \gamma} k_g < 4 \pi \, .
\end{equation}
\item[(C2)]
 For any $\epsilon > 0$,  if $R$ is sufficiently large, then
$\partial \Gamma_{\gamma , R} \setminus \gamma$ intersects $\{ x_3
< \epsilon \, R \}$.
\end{enumerate}
Since the statement is scale invariant, we can normalize so that
$\gamma \subset B_{1}$. We will take $R$ much larger than the
constant $\bar{R}$  given by Lemma \ref{l:ftc} so that \eqr{e:(D)}
holds for $R/2-1$.

The key point for proving (C1) is that the intrinsic one-sided
curvature estimate, Corollary \ref{t:one-sided}, gives a constant
 $\mu >0$ so that if $\epsilon < \mu$ and $y \in \{ |x_3| < \mu \, R \}
 \cap \partial \Gamma_{\gamma , R} \setminus \gamma$, then
\begin{equation}    \label{e:ceq}
    \sup_{\cB_{R/4}(y)}
    |A|^2 \leq C' \, \mu^2 \, R^{-2} \, .
\end{equation}
Note that to apply Corollary \ref{t:one-sided} here, we used Lemma
\ref{l:ftc} to see that $\cB_{R/2}(y)$ is a topological disk.
  The
claim (C1) follows easily from \eqr{e:ceq}.  Namely, first choose
a point $y_0 \in \{ x_3 < \epsilon \, R \} \cap \partial
\Gamma_{\gamma , R} \setminus \gamma$ and observe that the
curvature bound \eqr{e:ceq} allows us to apply the gradient
estimate to the positive harmonic function $x_3 + \epsilon \, R$
on $\cB_{R/4} (y_0)$ to get
\begin{equation}    \label{e:geq}
    \sup_{\cB_{R/8}(y_0)}
    |\nabla_{\Gamma} x_3| \leq C  \, \epsilon \, .
\end{equation}
The bound \eqr{e:geq} implies that the ball $\cB_{R/8}(y_0)$ is
graphical and moreover is contained in the slab $\{ |x_3| \leq C
\, \epsilon \, R \}$.  In particular, for $\epsilon
> 0$ sufficiently small, we can repeat this process to get a
chain of balls $\cB_{R/8}(y_i)$  with $y_i \in \partial
\Gamma_{\gamma , R} \cap \{ |x_3| < \mu \, R \}$ and so that
$\cup_i \cB_{R/8} (y_i)$ forms a graph which
 circles the $x_3$-axis. The intersection of this
 graph with the cylinder $\{ x_1^2 + x_2^2 = R^2 \}$
  contains a graph over the circle $\partial D_R$.  Since $\Gamma_{\gamma , 2R}$ is compact,
  the graph
  \begin{equation}
    \{ |x_3| < \mu \, R \} \cap \{ x_1^2 + x_2^2 = R^2 \} \cap
    \Gamma_{\gamma , 2R}
\end{equation}
 cannot spiral forever and,
hence, this graph closes up. Finally, the curvature bound
\eqr{e:ceq} and gradient bound for the graph imply a pointwise
bound for the geodesic curvature of $\partial \Gamma_{\gamma , R}$
and integrating this pointwise bound gives \eqr{e:cbd}.

To prove the second claim (C2), we will use catenoid barriers and
the strong maximum principle to argue by contradiction.  Suppose
therefore that $\epsilon > 0$ and
\begin{equation}    \label{e:stt}
    \partial \Gamma_{\gamma , R} \setminus \gamma  \subset  \{ x_3 > \epsilon \, R \} \, ,
\end{equation}
Let $\Cat$ denote the standard catenoid  ($\Cat = \{ \cosh^2 (x_3)
= x_1^2 + x_2^2 \} $) so that
\begin{equation}    \label{e:catslab}
    \{ x_1^2 + x_2^2 \leq 3 R \} \cap \Cat  \subset \{ |x_3| \leq
    \cosh^{-1} (3R) \} \, .
\end{equation}
Consider the one-parameter family of vertically translated
catenoids $\Cat_t = \Cat + (0,0,t)$ and observe that $\Cat_{-2R}
\cap \Gamma_{\gamma , R} = \emptyset$.  Furthermore, when $R$ is
large, \eqr{e:stt} and \eqr{e:catslab} imply that $\Cat_t \cap
\partial \Gamma_{\gamma , R} = \emptyset$ for every $t \leq 5 \, \cosh^{-1}
(3R)$.  Here we used \eqr{e:stt} to deal with the outer boundary
while the inner boundary $\gamma$ came for free since it is
contained in $B_1$.  By the strong maximum principle, there cannot
be a first $t \leq 5 \, \cosh^{-1} (3R)$ where $\Cat_t$ intersects
$\Gamma_{\gamma , R}$ and hence for $t \leq 5 \, \cosh^{-1} (3R)$
we have
\begin{equation}
    \Cat_{t} \cap \Gamma_{\gamma , R} = \emptyset \, .
\end{equation}
Arguing similarly give that a horizontal translation of $\Cat_{3
\, \cosh^{-1} (3R)}$ by a distance $2R$ cannot intersect
$\Gamma_{\gamma , R}$. However, this horizontally translated
catenoid separates $B_1$ and $\{ x_3 > \epsilon \, R \}$ in $B_R$
and hence separates the components of $\partial \Gamma_{\gamma ,
R}$, giving the desired contradiction.
\end{proof}

\subsection{The proof of Corollary \ref{c:finite}}

Both Corollary \ref{c:finite} and Corollary \ref{c:3rd} will use
the following weak chord arc property for annuli (cf. Proposition
\ref{t:2}):

\begin{Lem}     \label{t:2an}
Let $\Gamma$ be as in Lemma \ref{l:ftc}. There exist constants
$\tilde{R}$  and $\delta > 0$ so that for all intrinsic tubular
neighborhoods
 $T_{R}(\gamma)$ of $\gamma$ in $\Gamma$ with $R\geq \tilde{R}$,
 the component $\Gamma_{\gamma,\delta \, R}$ of
$B_{\delta \, R} \cap \Gamma$ containing $\gamma$ satisfies
 \begin{equation}   \label{e:t2an}
    \Gamma_{\gamma,\delta \, R} \subset T_{R/2} (\gamma) \, .
 \end{equation}
 Here $\tilde{R}$ depends on
$\Gamma$ but $\delta$ does not.
\end{Lem}

\begin{proof}
Let $\bar{R}$ be the constant given by   Lemma \ref{l:ftc}  so
that \eqr{e:(D)} holds.  We can now directly follow the proof of
claim \eqr{e:claim} in the proof of Proposition \ref{p:key}
 to get \eqr{e:t2an}.  This  requires one modification to get that intrinsic
 subballs are weakly chord arc. Namely,
 rather than using condition (A') there, we  use \eqr{e:(D)} to first
 see that the intrinsic subballs
 are disks and then apply Proposition \ref{t:2} to these disks.
\end{proof}

The weak chord arc property given by Lemma \ref{t:2an} implies the
necessary compactness needed to apply Lemma \ref{l:annuli1} and
gives that embedded minimal annuli in a halfspace have finite
total curvature:

\begin{Cor}     \label{c:annuli}
Let $\Gamma$ be as in Lemma \ref{l:ftc}.  If $\Gamma$ is contained
in a halfspace, then $\Gamma$ has finite total curvature.  It
follows that $\Gamma$ is asymptotic to a plane or half of a
catenoid.
\end{Cor}

\begin{proof}
Lemma \ref{t:2an} implies that, for every $R$, the component
$\Gamma_{\gamma , 2R}$ of $B_{2R} \cap \Gamma$ containing $\gamma$
is compact.  Hence, we can apply Lemma \ref{l:annuli1} to
$\Gamma_{\gamma , 2R}$ for $R$ sufficiently large to get
\begin{equation}    \label{e:lrti}
     \int_{\Gamma_{\gamma , R}} |A|^2 \leq 2 \, \int_{\gamma}
    k_g + 8\pi \, .
\end{equation}
 As $R$ goes to infinity, the $\Gamma_{\gamma , R}$'s exhaust $\Gamma$ and
 hence \eqr{e:lrti} bounds the total curvature of $\Gamma$.  The
 second statement follows since the annulus $\Gamma$ is also embedded
 (see, e.g., \cite{Sc2}).
\end{proof}

 Corollary \ref{c:finite}, and hence Calabi's conjectures
for surfaces with finite topology,  now follow easily from
Corollary \ref{c:annuli}:

\begin{proof}
(of Corollary \ref{c:finite}). Observe first that an embedded
minimal surface $\Sigma$ with finite topology in a halfspace has
finite total curvature. This is because such a  $\Sigma$ can be
written as the union of a compact piece $\Sigma_0$ which may have
nonzero genus and a finite collection of non-compact annuli
$\Gamma_1 , \dots , \Gamma_k$ each of which contains one of its
boundary components. Clearly, $\Sigma_0$ has finite total
curvature since it is compact. Furthermore, each $\Gamma_i$ has
finite total curvature by Corollary \ref{c:annuli}, so we conclude
that $\Sigma$ itself has finite total curvature.

Finally, since $\Sigma$ has finite total curvature, \cite{Hu}
implies that $\Sigma$ is parabolic (in the sense that any positive
harmonic function is constant).  Therefore the positive harmonic
function $x_3$ is constant on $\Sigma$ and $\Sigma$ must be a
plane as claimed.
\end{proof}

\subsection{The proof of Corollary  \ref{c:3rd}: Properness}

The properness of embedded minimal surfaces with finite topology
will be an almost immediate consequence of properness of embedded
annuli that we will show next.  As in the case of disks, the weak
chord arc property given by Lemma \ref{t:2an} applies only to one
component and therefore does not directly give properness.

\begin{Pro}         \label{p:pan}
Let $\Gamma$ be as in Lemma \ref{l:ftc}.  Then $\Gamma$ must be
proper.
\end{Pro}

\begin{proof}
 The proposition follows from the
following claim: For every radius $R > 0$,  there is a constant
$S_R > R$ (depending on both $R$ and $\Gamma$) so that
\begin{equation}    \label{e:claimsr}
  B_R \cap \Gamma \subset \Gamma_{\gamma , S_R} \, .
\end{equation}
Here, as in Lemma \ref{t:2an}, $\Gamma_{\gamma , S_R}$ denotes the
component of $B_{S_R} \cap \Gamma$ containing $\gamma$.  To get
the proposition from \eqr{e:claimsr}, simply apply Lemma
\ref{t:2an} (for $R$ large) to get
\begin{equation}   \label{e:t2an3}
    \Gamma_{\gamma, S_R} \subset T_{S_R/(2\delta)} (\gamma) \, ,
\end{equation}
and  observe that the (closure of the) intrinsic tubular
neighborhood $T_{S_R/(2\delta)} (\gamma)$ is compact.

The rest of the proof is to establish \eqr{e:claimsr}.  We will do
this by contradiction; suppose therefore that $R > 0$ is fixed,
$\gamma \subset B_R$,  and $y_i$ is a sequence of points in $B_R
\cap \Gamma$ with
\begin{equation}        \label{e:notyi}
    y_i \notin \Gamma_{\gamma , i \, R} \, .
\end{equation}
We will show that \eqr{e:notyi} implies that $\Gamma$ has finite
total curvature and then get a contradiction from this.

The first step is to find large graphical regions in $\Gamma$.
Observe that, by the triangle inequality,
\begin{equation}    \label{e:diyi}
    d_i = \dist_{\Gamma} (y_i , \gamma) \geq (i-1) \, R \, .
\end{equation}
Since $i \to \infty$, it follows from \eqr{e:diyi}  that for any
$J$ we can choose indices $i_1$ and $i_2$ so that
\begin{align}    \label{e:theJ1}
    d_{i_1} > 2 J {\text{ and }} d_{i_2} &> 2 J \, ,  \\
     \dist_{\Gamma}
    (y_{i_1} , y_{i_2} ) &> 2 J \, . \label{e:theJ2}
\end{align}
When $J$ is large, Lemma \ref{l:ftc} and \eqr{e:theJ1} imply that
the intrinsic balls $\cB_{J}(y_{i_1})$ and  $\cB_{J}(y_{i_2})$ are
topological disks; and disjoint by \eqr{e:theJ2}.  The  one-sided
curvature estimate now implies that
 $\cB_{J}(y_{i_1})$ and $\cB_{J}(y_{i_2})$   contains
a graph $\Gamma_1$ and $\Gamma_2$, respectively,  over a disk of
radius $c \, J$ with small gradient $\leq \tau$  and $y_{i_j} \in
\Gamma_j$ (where $c$ depends on $\tau$).  To prove this, first
apply Proposition \ref{t:2} to see that the intrinsic balls are
weakly chord arc and then apply Corollary \ref{c:barrier} to get a
curvature bound.

The second step is to use the large graphical region  to show that
$\Gamma$ has finite total curvature. Namely, for $J$ large, Lemma
\ref{t:2an} implies that the component $\Gamma_{\gamma , cJ}$ of
$B_{cJ} \cap \Gamma$ containing $\gamma$ is compact. Moreover,
since $\Gamma$ is embedded, the graph $\Gamma_1$ forces
$\Gamma_{\gamma , cJ}$ to be contained in a halfspace
\begin{equation}
    \Gamma_{\gamma , cJ} \subset \{ x_3 > -R - c \, \tau \, J \}
    \, .
\end{equation}
(Here we have assumed that $\Gamma_1$ is beneath $\gamma$; this
can be arranged after possibly reflecting across $\{ x_3 = 0 \}$.)
 For $\tau > 0$ small,
we can apply Lemma \ref{l:annuli1} to get a  bound for the total
curvature of $\Gamma_{\gamma , c J/2}$  which is independent of
$J$. It follows that $\Gamma$ has finite total curvature since the
$\Gamma_{\gamma , c J/2}$'s exhaust $\Gamma$ as $J \to \infty$.

Finally,  as in the
 proof of Lemma \ref{l:ftc}, we conclude that $\Gamma$ is
 asymptotic to either a plane or half of a catenoid since it has
 finite total curvature.  However, in either case, \eqr{e:claimsr}
 clearly holds.  This contradiction establishes the claim
 \eqr{e:claimsr} and thus completes the proof.
\end{proof}

The properness of embedded minimal surfaces with finite topology
now follows easily:

\begin{proof}
(of Corollary \ref{c:3rd}).  Write the embedded minimal surface
with finite topology $\Sigma$ as the union of a compact piece
$\Sigma_0$  and a finite collection of non-compact annuli
$\Gamma_1 , \dots , \Gamma_k$ each of which contains one of its
boundary components.  Proposition \ref{p:pan} implies that each
annulus $\Gamma_i$ is proper and hence so is $\Sigma$.
\end{proof}

\appendix

\section{Multi-valued graphs} \label{s:a1}

To make the notion of multi-valued graph precise,  let $\cP$ be
the universal cover of the punctured plane $\CC\setminus \{0\}$
with global polar coordinates $(\rho, \theta)$ so $\rho>0$ and
$\theta\in \RR$.  An {\it $N$-valued graph} on the annulus
$D_s\setminus D_r$ is a single valued graph of a function $u$ over
\begin{equation}
    \{(\rho,\theta)\,|\,r< \rho\leq s\, ,\, |\theta|\leq
    N\,\pi\} \, .
\end{equation}
 For working purposes, we generally
think of the intuitive picture of a multi-sheeted surface in
$\RR^3$, and we identify the single-valued graph over the
universal cover with its multi-valued image in $\RR^3$.

The multi-valued graphs that we   consider in this paper will all
be embedded, which corresponds to a nonvanishing separation
between the sheets. Here the {\it separation} is the function
\begin{equation}
w(\rho,\theta)=u(\rho,\theta+2\pi)-u(\rho,\theta)\, .
\end{equation}
If $\Sigma$ is the helicoid [i.e., $\Sigma$ can be parametrized by
$(s\,\cos t,s\,\sin t,t)$ where $s,\,t\in \RR$], then
$\Sigma\setminus \{x_3-\text{axis}\}=\Sigma_1\cup \Sigma_2$,
 where $\Sigma_1$, $\Sigma_2$ are $\infty$-valued graphs on
$\CC\setminus \{0\}$. $\Sigma_1$ is the graph of the function
$u_1(\rho,\theta)=\theta$ and $\Sigma_2$ is the graph of the
function $u_2(\rho,\theta)=\theta+\pi$.  ($\Sigma_1$ is the subset
where $s>0$ and $\Sigma_2$ the subset where $s<0$.)  In either
case the separation $w=2\,\pi$.

\section{The proof of Lemma \ref{l:ca}}     \label{s:a4}

We will next include the proof of Lemma \ref{l:ca}.  This lemma is
modelled on lemma II.2.1 in \cite{CM6}.  The proof follows that of
lemma II.2.1 in \cite{CM6} with very minor changes, but we include
it here for completeness.

\begin{proof}
 (of Lemma \ref{l:ca}).
Let $C_0 > 2$ be given by Lemma \ref{l:halfstable}.  We will show
that there exists $n$ depending on $C_a$ so that
 \begin{equation}        \label{e:claim2}
        \Sigma_{0, \frac{R}{C_0} } \subset
        \cB_{n \, R}(0) \, .
\end{equation}
 To prove this, we will argue by contradiction; so suppose that
 \eqr{e:claim2}
fails for some large $n$.  Consequently, we get a curve
\begin{equation}
    \sigma \subset \Sigma_{0, \frac{R}{ C_0} } \subset B_{\frac{R}{
    C_0}}(0)
\end{equation}
 from $0$ to a point in $\partial \cB_{n \, R}(0)$. For $i=1, \dots ,
n$, fix points
\begin{equation}
    z_i \in \partial \cB_{i \,  R}(0) \cap \sigma \, .
\end{equation}
  It follows that the intrinsic balls $\cB_{R/2}(z_i)$:
\begin{itemize} \item
 Are
disjoint.
\item Have
centers in $B_{\frac{R}{ C_0}}(0)$.
\item Do not intersect $\partial \Sigma$.
\end{itemize}
Since the $n$ points $\{ z_i \}$ are all in the Euclidean ball
$B_{\frac{R}{
    C_0}}(0) \subset \RR^3$,
    there exist integers $i_1$ and $ i_2$ with
\begin{equation}    \label{e:4272}
    0< | z_{i_1} - z_{i_2}| < C \,  n^{-1/3} \, R \, .
\end{equation}
 Note that  \eqr{e:cacb} gives a uniform curvature bound on the
balls $\cB_{R/2}(z_{i_1})$ and $\cB_{R/2}(z_{i_2})$.  Therefore,
 Lemma \ref{l:halfstable} implies that, for $n$
sufficiently large (so the centers $z_{i_1}$ and $z_{i_2}$ are
extrinsically close), then we get for $j=1$, $2$ that
\begin{equation}    \label{e:this}
    B_{ \frac{R}{C_0} }(0) \cap
\partial \cB_{ \frac{11 R}{2C_0} }(z_{i_j})  = \emptyset \, .
\end{equation}
(Here we used that $B_{R/C_0}(0)  \subset B_{5R/C_0}(z_{i_j})$
because $z_{i_j} \in B_{R/C_0}(0)$.)
 Since the curve $\sigma$
must intersect $\partial \cB_{11 R/ (2C_0)}(z_{i_j})$,
\eqr{e:this}  contradicts
 the fact that the curve $\sigma$ is contained in
the ball $B_{  R / C_0 }(0)$.  This  contradiction proves
\eqr{e:claim2} and consequently
  gives the lemma.
\end{proof}

\end{document}